\documentclass[12pt,a4paper]{article}

\usepackage{makeidx}         
\usepackage{graphicx}        
\usepackage{psfrag}                             
\usepackage{multicol}        
\usepackage{amsmath}
\usepackage{amssymb}
\usepackage{amsfonts}
\usepackage{comment}
\usepackage{latexsym}
\usepackage[T1]{fontenc}
\usepackage{subfigure}                      
\usepackage{float}
\usepackage{mathrsfs}
\usepackage{mathptmx}
\usepackage{upref}      
\usepackage{hyperref}
\usepackage{color}
\usepackage{amsthm}
\usepackage{url}

\newtheorem{theorem}{Theorem}[section]

\newtheorem{lemma}[theorem]{Lemma}

\let\eps\varepsilon

\newcommand{\ee}{\varepsilon}

\begin{document}
\title{\bf Asymptotics for the spectrum of the Laplacian  in thin bars with varying cross sections}

\author{Pablo Benavent-Ocejo\footnote{Universidad de Cantabria, Av. Los Castros, 39005, Santander, Spain, \url{pablo.benavent@alumnos.unican.es}}, \, Delfina G\'omez\footnote{Departamento Matem\'aticas, Estad\'istica y Computaci\'on, Universidad de Cantabria, Av. Los Castros, 39005, Santander, Spain, \url{gomezdel@unican.es}} \, Maria-Eugenia P\'erez-Mart\'inez\footnote{ Departamento de Matem\'atica Aplicada y Ciencias de la Computaci\'on, Universidad de Cantabria, Av. Los Castros, 39005, Santander, Spain, \url{meperez@unican.es}}\, \footnote{Corresponding author}}
\date{}

\maketitle

\noindent
\abstract{
\noindent
We  consider spectral problems for Laplace operator in  3D rod structures with a  small cross section of diameter $O(\varepsilon)$, $\varepsilon$ being a positive parameter. The boundary conditions are Dirichlet  (Neumann, respectively)  on  the bases of this structure and Neumann on the lateral boundary.   As   $\varepsilon\to 0$, we show the convergence of the spectrum with conservation of the multiplicity towards that of a 1D spectral  model  with Dirichlet (Neumann, respectively) boundary conditions.  This 1D model may arise in diffusion or vibrations models of nonhomogeneous media with different physical characteristics and it takes into account the geometry of the 3D domain.  We deal with the low frequencies and  the approach to eigenfunctions  in the suitable Sobolev spaces is also outlined.  }


\section{Introduction}\label{sec:ben1}

In this paper, we address the asymptotic behavior of the eigenvalues and eigenfunctions  for the Laplace operator in  thin rod structures,  when the diameter of the {\em cross / transverse} section  tends to zero. Namely, 3D  domains of size $O(1)$ along the longitudinal direction and  $O(\varepsilon)$ in the two other directions which are referred to as {\em transverse directions} (cf. Figures~\ref{Fig_domain}-\ref{figura_boveda}). The boundary conditions are Neumann on the lateral surface while they can be Dirichlet or Neumann on the parallel  bases (cf. \eqref{formulacionuno_b}-\eqref{fN}). The spectrum of the problem being discrete, we establish the order of magnitude of the so-called {\em low frequencies} and show that the limit problem is a 1D model which takes into account the geometrical characteristics of  the rod and the boundary conditions, while the eigenfunctions become strongly  oscillating functions in the longitudinal direction as the eigenvalue number increases and $\varepsilon\to 0$ (cf. Figures \ref{figura_mixed}-\ref{figura_mixed2}).

\medskip
\noindent
Rod    structures appear in many engineering constructions or engineering devices containing  bar or thin tube structures, but  from the  mathematical viewpoint there are  gaps in the   description of eigenvalues and eigenfunctions for 3D structures, in their dependence of the small parameter; specially the models with the boundary conditions that we address here are   open problems in the literature of Applied Mathematics. The interest from the dynamical viewpoint is evident, both for models arising in diffusion  or vibrations of tube structures and / or  multistructures (cf. \cite{AmosovGomezPanasenkoPerez2024,Pa_Pi_book,Nazarov2024,Nazarov2025}).

\medskip
\noindent
Explicit computations on  prism-like  structures have been performed in \cite{articulo} and show the different asymptotic  behavior of the eigenvalues and eigenfunctions depending on the boundary conditions.  They are important to  enlighten  the order of magnitude of  the low frequencies, its asymptotic behavior, as $\varepsilon\to 0$, and the behavior of the associated eigenfunctions. These  computations also reveal that in order to capture oscillations of the eigenfunctions different from the longitudinal ones   we need to deal with the high frequencies. Also they show how different the above-mentioned behavior is depending on the boundary conditions, see Remark \ref{Darticulo} in this connection.  Additionally, when explicit computations do not work due to the geometry of the domain, experimentally, numerical methods  using finite elements (cf. the PDE Toolbox of Matlab) show this different behavior in     \cite{articulo} but also allow  us to observe the  numerical instabilities that make  it impossible to approximate   the eigenvalues and eigenfunctions as $\varepsilon$ becomes smaller.

\medskip
\noindent
This is why the asymptotic study that we perform here becomes necessary. As a matter of fact, for mixed boundary conditions (Dirichlet on the parallel bases), the order of magnitude of the low frequencies and the limit problem has been obtained in \cite{articulo} using the minimax principle and taking limits in the variational formulation using special test functions. The convergence with conservation of the multiplicity was left as an open problem that we address here
(cf. \eqref{conservation}). Also, we show the convergence of the eigenvalues and the associated eigenfunctions with conservation of the multiplicity for  the Neumann problem (cf. \eqref{conservation2}).
We gather the main convergence results in Theorems~\ref{TheoremM} and \ref{TheoremN}.
In order to do it we use a  result from the spectral perturbation theory (cf. Lemma \ref{LemaOl}) which is well adapted to obtain  spectral convergence when operators and Hilbert spaces depend on  the small parameter. It allows us to show the convergence of the eigenvalues and eigenfunctions when the convergence of the solutions of the  associated stationary problems is obtained. Thus, we focus on stationary problems and introduce   techniques based on reduction of dimensions and weighted Sobolev spaces to obtain the limit problems.

\medskip
\noindent
Let us recall that, up to our knowledge, in the framework of thin bars/rods for the dimension 3 of the domain,  the spectral mixed boundary value problem had not been considered previously in the literature, while the Dirichlet and Neumann ones  had  but for different  geometrical configurations and using different techniques that we briefly explain.

\medskip
\noindent
For a fixed cross section, and different geometries, Dirichlet boundary conditions have been considered in \cite{CardoneDranteNazarov2010} and  \cite{BorisovCardone2011} for 3D rods
by means of different techniques.  Mixed boundary conditions (but the Dirichlet one  being imposed on the lateral boundary) and Neumann conditions  have been  addressed in \cite{CardoneDranteNazarov2010} focussing on the localization phenomena for the eigenfunctions near some perturbed ends of the rod while the cross  section stays constant.    See  \cite{FriedlanderSolomyak}  for other localization phenomena in 2D  thin domains, and    \cite{BorisovCardone2011} and \cite{CardoneDranteNazarov2010} for further references.       Spectral problems for 3D tube structures with Neumann boundary conditions have been considered in \cite{AmosovGomezPanasenkoPerez2023,AmosovGomezPanasenkoPerez2024} where  the asymptotic partial decomposition of the domain method   is used to approach  the problem by a spectral problem in a domain   with hybrid dimensions (3D and 1D); the tubular part of the structure somehow replaced by the 1D structure has a  constant cross section. Also, \cite{PanasenkoPerez2007} uses this technique of domain decomposition for a mixed boundary value problem in a thin planar domain. For thin 2D rod structures and 3D like films structures with only one of the dimensions smaller than the other or an oscillating boundary,   we refer to \cite{Nazarov2013,NazarovPerezTaskinen2016, ArrietaNakasatoVillanueva2025}  and references therein. The junction of rod structures or  thin films has been addressed in \cite{GaudielloSili2007, GaudielloGomezPerez2023, Chesnel}.

\medskip
\noindent
Finally, let us describe the structure of the paper. Section \ref{sec:ben2} contains the statement of the problem under consideration and the limit problems.
Some background and preliminary results   useful for proofs  are in Section  \ref{sec:ben2.1}. Section  \ref{sec:ben3} addresses the proof of the convergence in the case of mixed boundary conditions while   Section  \ref{sec:ben4}  contains the proof of the convergence for the Neumann Laplacian.

\section{Setting of the problem}\label{sec:ben2}
For each $\ee \in (0,1)$,  let us introduce  $G_\ee$ a rod-type domain   ({\em rod} in short), which for simplicity we can assume  located along the $x_1$- axis as follows.

\noindent
First, we consider an open bounded     domain  $G\subset\mathbb{R}^3$ with a Lipschitz boundary  that admits the following representation
\begin{equation}
    G=\bigcup_{x_1\in(l_0,l_1)}\{(x_1,x_2,x_3): (x_2,x_3)\in D_{x_1}\},
    \label{dominioG}
\end{equation}
where $l_0,l_1\in \mathbb{R}$,  $l_0\leq 0< l_1$ and   $D_{x_1}$ denotes the cross section at the point $x_1\in (l_0,l_1)$;   namely, for any fixed $x_1\in(l_0,l_1)$,  $D_{x_1}$ is an open domain of the   $x_2x_3$-plane depending on $x_1$ and such that $(0,0)\in \overline{D_{x_1}}$. Moreover, we assume that the function $\vert D_{x_1}\vert$  defining the areas of the cross sections of the domain $G$  is a piecewise continuous function in $[l_0,l_1]$  and satisfy
\begin{equation}
    0<c_0<\vert D_{x_1}\vert\leq c_1,\quad \forall x_1\in[l_0,l_1],
    \label{acotacionarea}
\end{equation}

\medskip
\noindent
We denote by $\partial G$ the boundary of $G$, which consists of the union of two flat faces, {\em the parallel  bases}, denoted by $\Gamma_0$ and $\Gamma_1$, and a lateral surface $\Gamma_l\subset\mathbb{R}^3$. Further specifying:   $$\overline{\Gamma}_0=\overline{G}\cap\{x_1=l_0\},\quad   \overline{\Gamma}_1=\overline{G}\cap\{x_1=l_1\}\quad \text{and} \quad \partial G=\overline{\Gamma}_0\cup\overline{\Gamma}_1\cup\overline{\Gamma}_l.$$ For simplicity, in the description it has been assumed that the origin of coordinates is in $\overline{G}$.

\begin{figure}
\centering
\includegraphics[width=4.5cm]{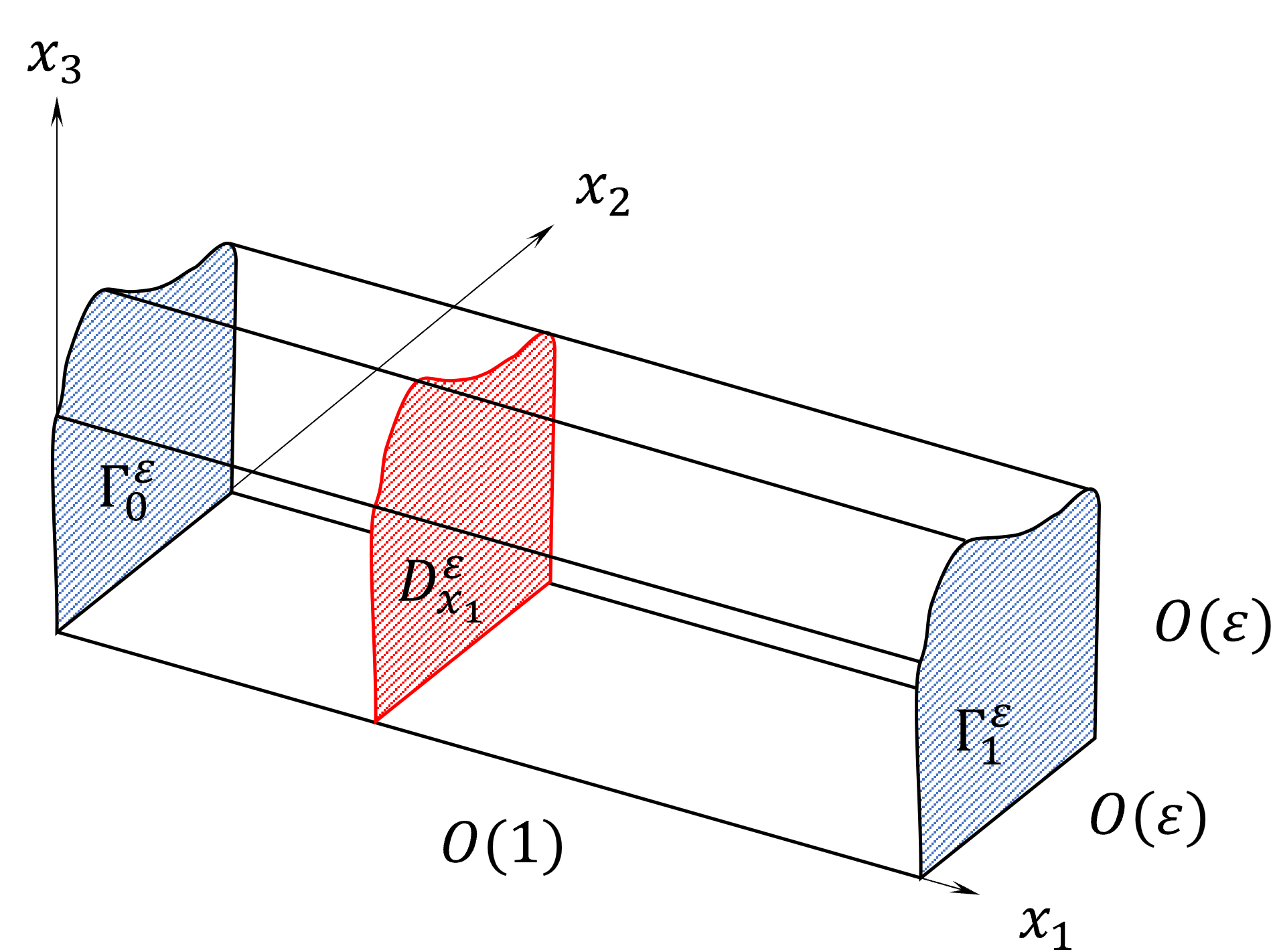}
\caption{\small Geometrical configuration of a  domain $G_\varepsilon$.}
\label{Fig_domain}
\end{figure}

\medskip
\noindent
Now, let $\varepsilon\in(0,1)$ denote a small parameter that we shall make to go to  $0$, and   consider
\begin{equation}\label{dominioGeps} G_\varepsilon:=\{(x_1,x_2,x_3):\Big(x_1,\frac{x_2}{\varepsilon},\frac{x_3}{\varepsilon}\Big)\in G\},\end{equation}
that is,
$$G_\varepsilon=\bigcup_{x_1\in(l_0,l_1)}\{(x_1,x_2,x_3):(x_2,x_3)\in D_{x_1}^\varepsilon\}\quad\text{ with }D_{x_1}^\varepsilon=\varepsilon D_{x_1}$$
We denote by $\Gamma_\varepsilon^D$ the two faces perpendicular to the $x_1$-axis, that is:
$$\Gamma_\varepsilon^D=\Gamma_0^\varepsilon\cup\Gamma_1^\varepsilon\quad \text{ with }\quad \overline{\Gamma_0^\varepsilon}=\overline{G}_\varepsilon\cap\{x_1=l_0\}\quad \text{ and }\quad \overline{\Gamma_1^\varepsilon}=\overline{G}_\varepsilon\cap\{x_1=l_1\},$$
and the lateral surface is $\Gamma_\varepsilon^l=\partial G_\varepsilon\setminus\overline{\Gamma_\varepsilon^D}$ (see Figure~\ref{Fig_domain}).

\medskip
\noindent
In $G_\ee$, we consider the spectral problem for the Laplacian consisting of the set of equations
\begin{eqnarray}
&\displaystyle -\Delta u_\varepsilon=\lambda_\varepsilon u_\varepsilon & \text{in }G_\varepsilon  \label{formulacionuno_a},\\   [0.4em]
&\displaystyle      \dfrac{\partial u_\varepsilon}{\partial n}=0 & \text{on }  {\Gamma_\varepsilon^l},  \label{formulacionuno_b}
\end{eqnarray}
and either the Dirichlet   conditions
\begin{equation}\label{fD} u_\varepsilon=0  \text{ on }\Gamma_\varepsilon^D,\end{equation}
or the Neumann  conditions
\begin{equation}\label{fN}  \dfrac{\partial u_\varepsilon}{\partial n}=0    \text{ on }\Gamma_\varepsilon^D.\end{equation}
\noindent
Above   $n$ stands for the outward unit normal vector to $\partial G_\varepsilon$, and  the pair $(\lambda_\varepsilon, u_\varepsilon)$ represents the eigenvalue  and the corresponding eigenfunction. In what follows \eqref{formulacionuno_a},\eqref{formulacionuno_b},\eqref{fD} is referred to as the mixed boundary value problem while \eqref{formulacionuno_a},\eqref{formulacionuno_b},\eqref{fN}  is   the Neumann problem (cf. also \eqref{formulacionproblemaNeumann}).

\medskip
\noindent
The limit eigenvalue problem  is defined by the set of equations
\begin{equation}\label{limitEDP}
         -\partial_{x_1}(|D_{x_1}| \partial_{x_1} U_0) = \lambda_0 |D_{x_1}| U_0, \quad x_1 \in (l_0,l_1), \end{equation}
         and either the Dirichlet condition (related to \eqref{fD})  \begin{equation}\label{limitD}
         U_0 (l_0) = 0, \quad U_0 (l_1) = 0, \end{equation}  or the Neumann one (related to \eqref{fN}) \begin{equation}\label{limitN}
         U_0' (l_0) = 0, \quad U_0'(l_1) = 0,  \end{equation}
       while $(\lambda_0, U_0)$ stands for  the spectral parameter  and the associated eigenfunction.

\begin{figure}[ht]
\centering
\includegraphics[width=4.5cm]{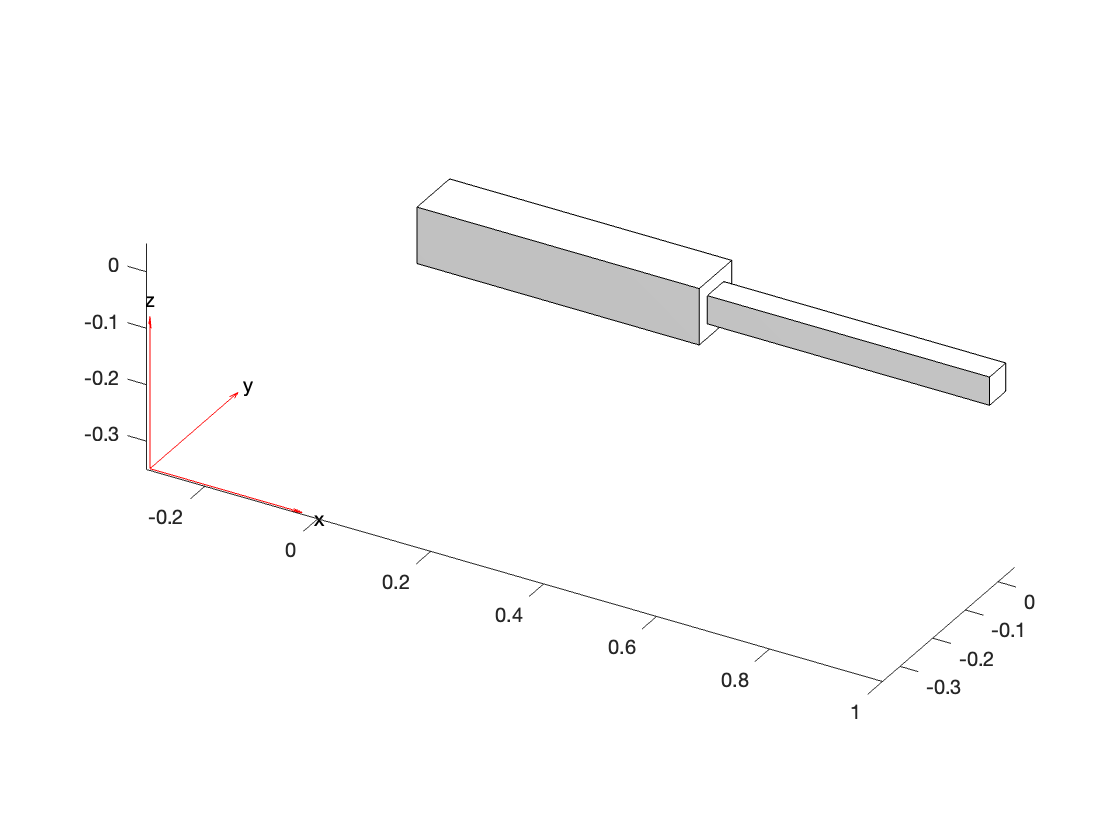}
\quad
\includegraphics[width=4.5cm]{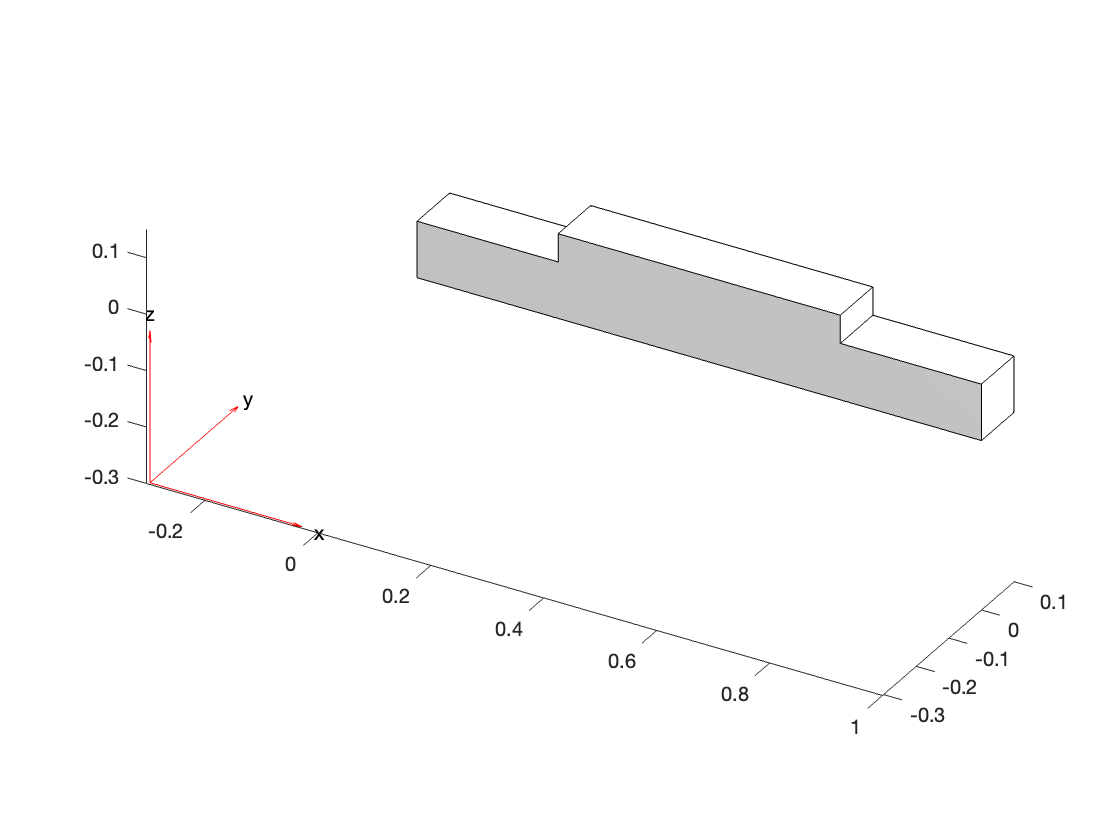}
\quad
\includegraphics[width=4.5cm]{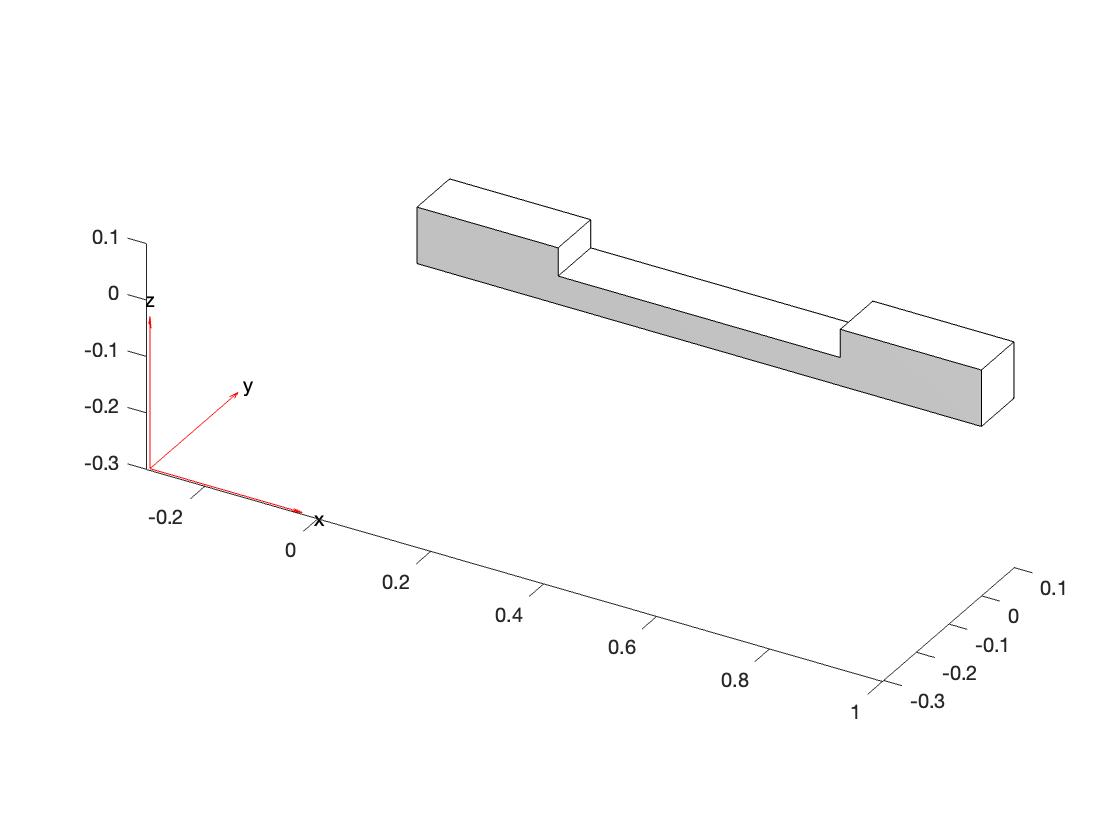}
\caption{\small Polyhedral domains which are the union of two cuboids  of different cross sections (of size $O(\ee)$ in both transverse directions)  or  union of   cuboids with different heights $O(\ee)$.}
\label{ej_unionprismas}
\end{figure}

\begin{figure}[ht]
\centering
\includegraphics[width=4.5cm]{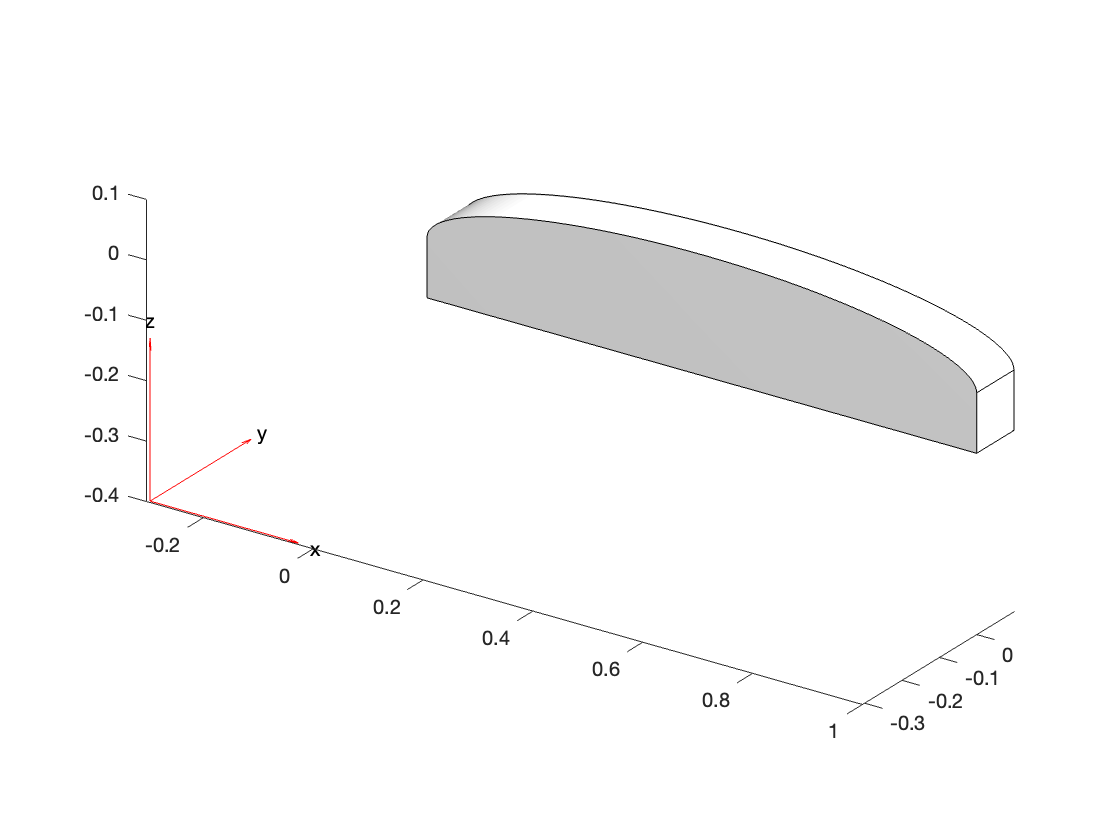}
\quad
\includegraphics[width=4.5cm]{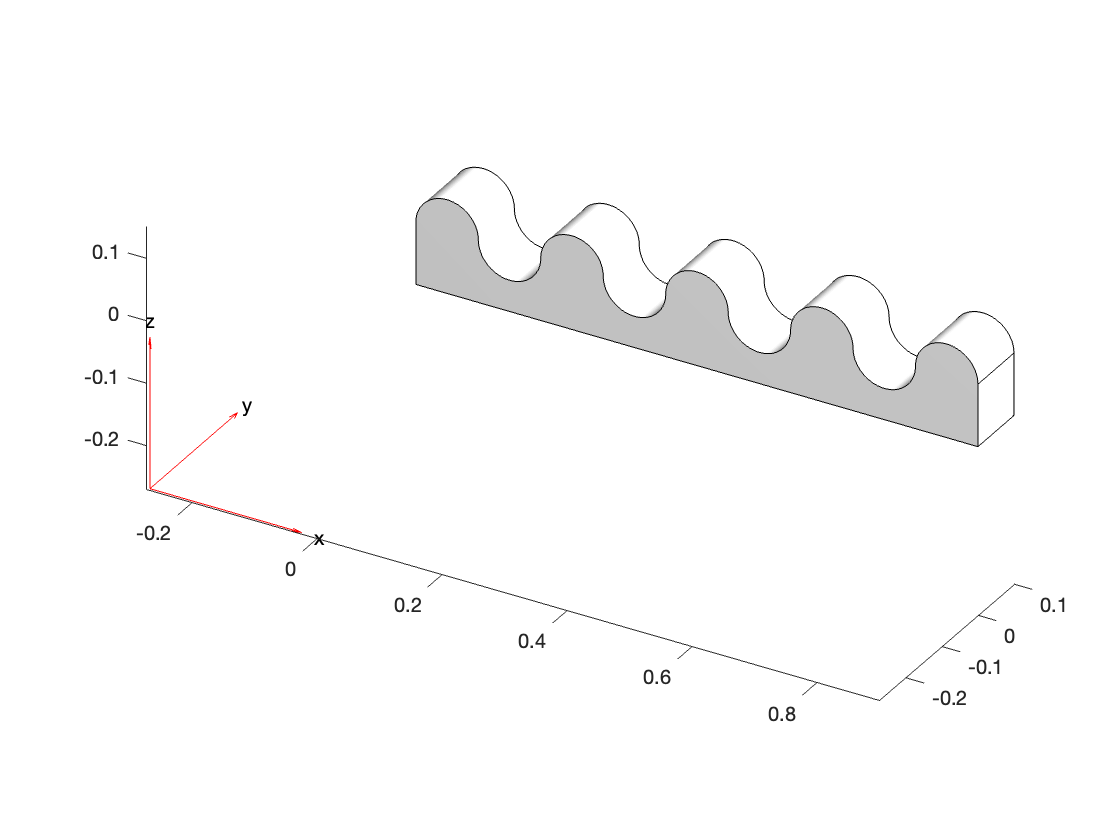}
\quad
\includegraphics[width=4.5cm]{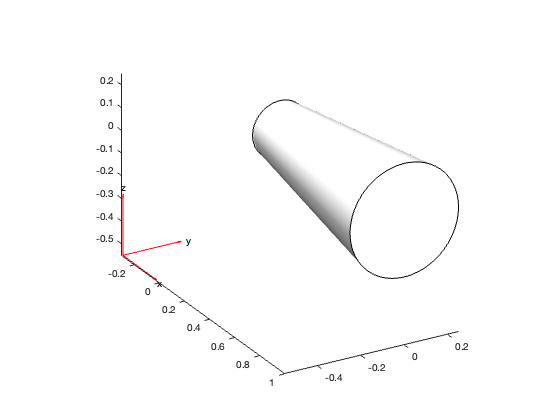}
 \caption{\small Dome-like top    or wavy-like top domain   of size $O(\ee)$ in both transverse directions and trumpet-like domain of diameter $O(\ee)$.}
 \label{figura_boveda}
\end{figure}

\begin{figure}[ht]
\begin{center}
\scalebox{0.15}{\includegraphics{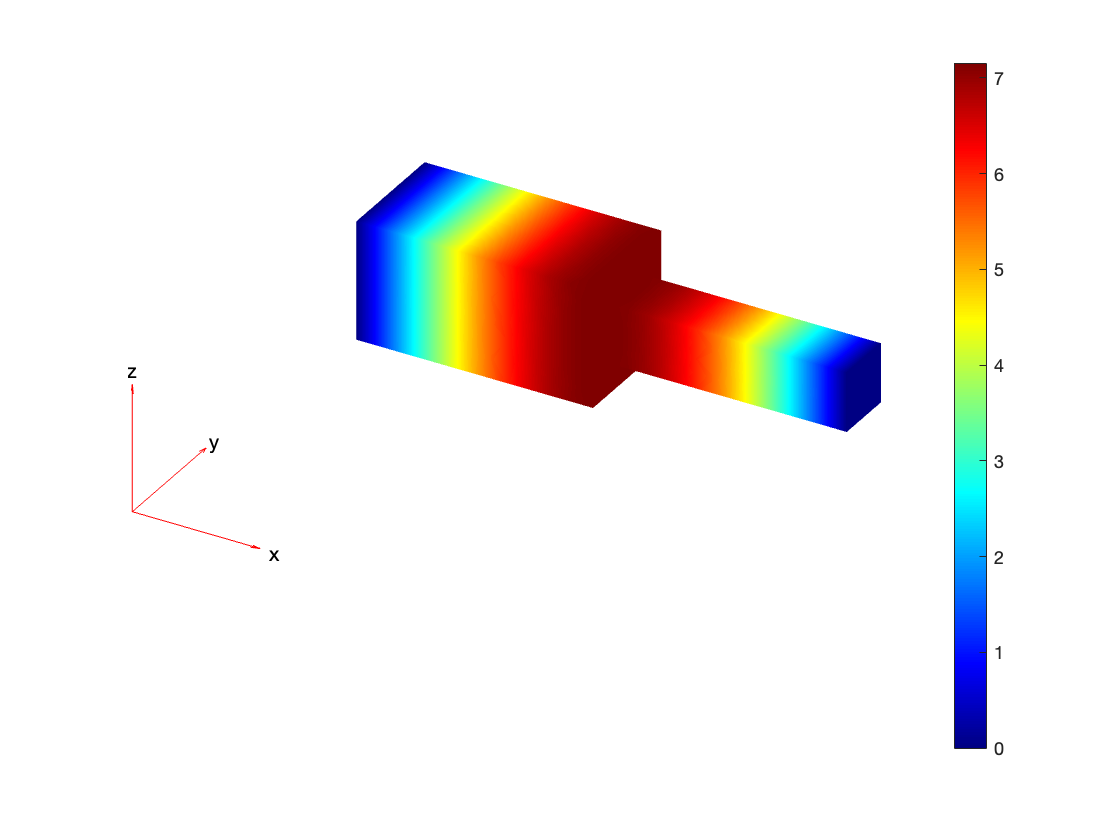}}
\scalebox{0.15}{\includegraphics{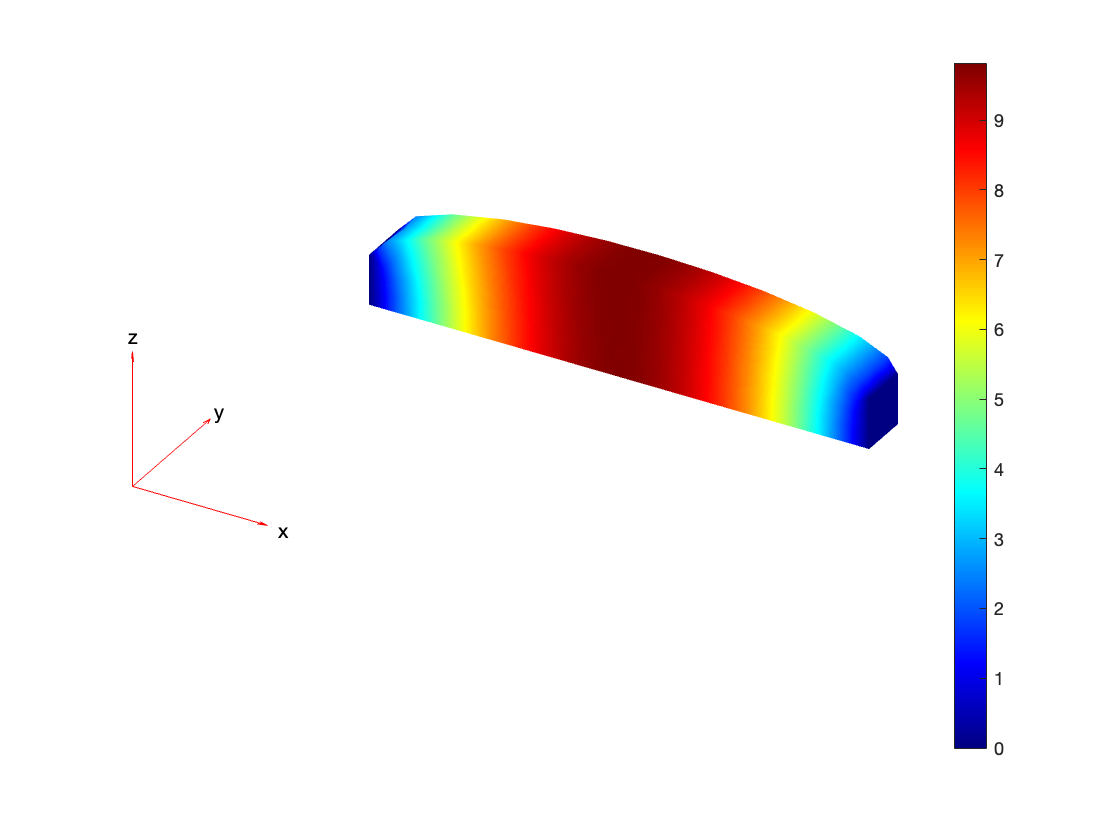}}
\scalebox{0.15}{\includegraphics{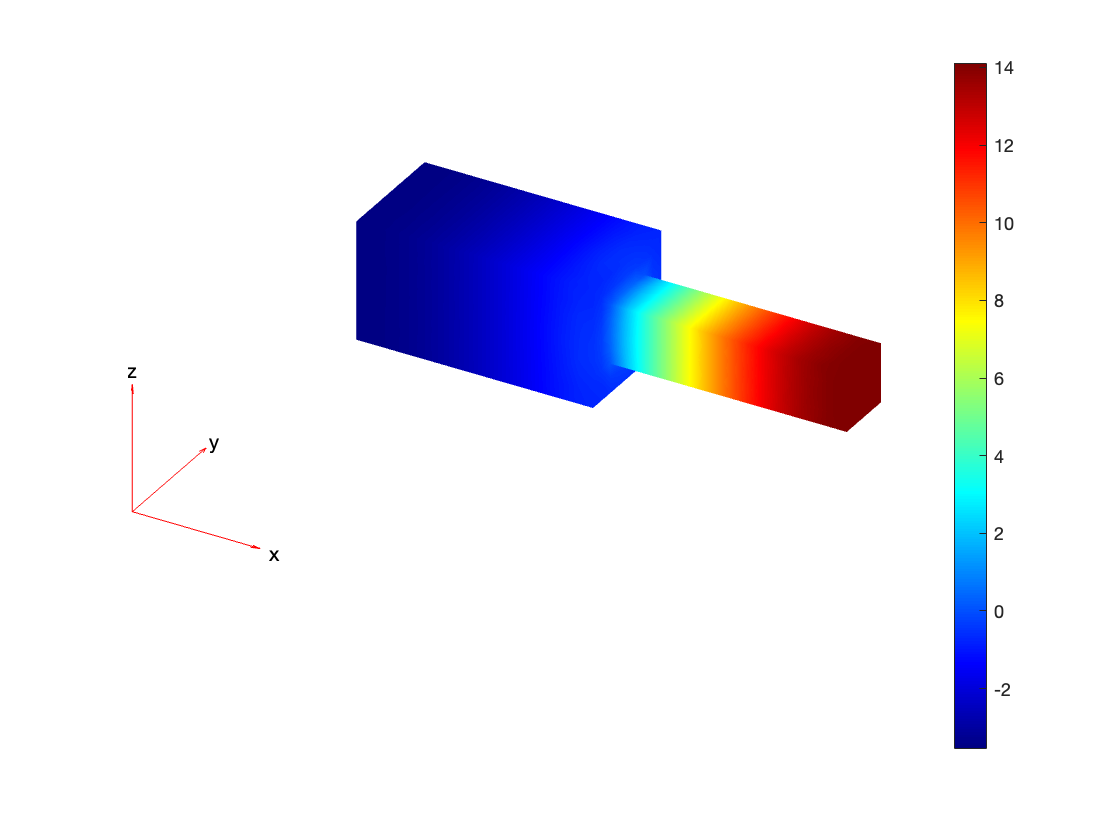}}
\scalebox{0.15}{\includegraphics{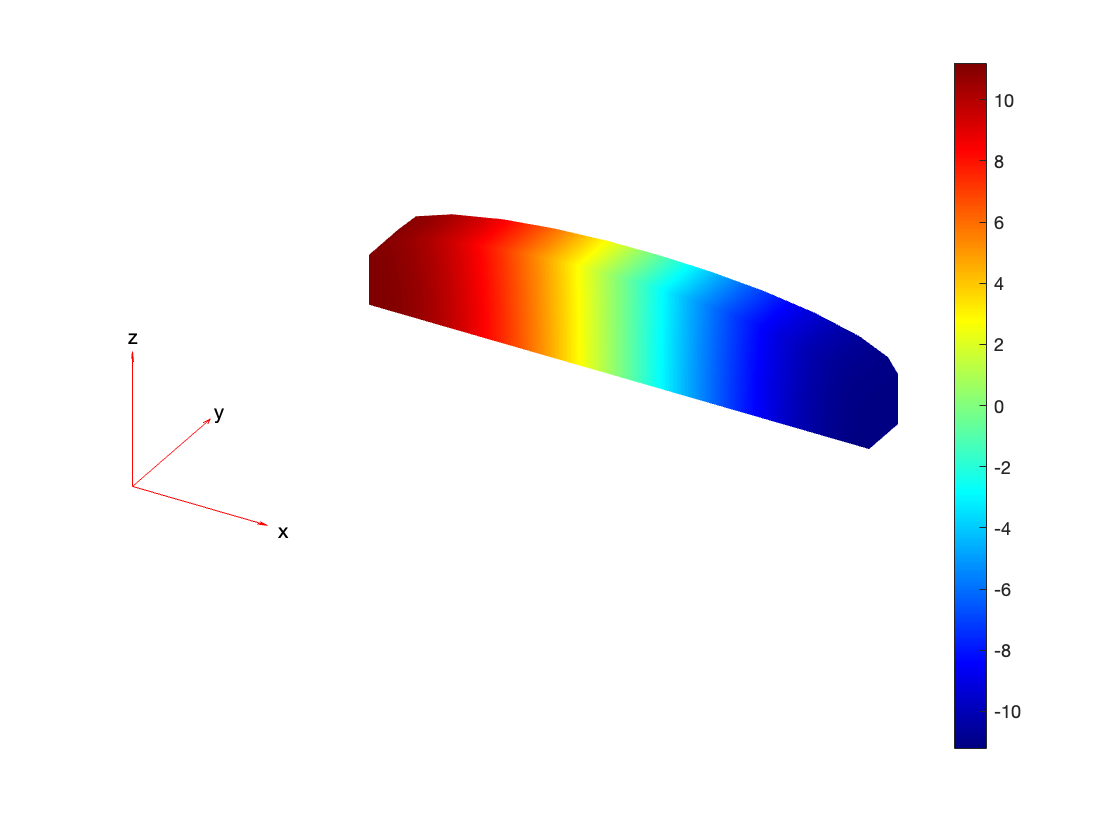}}
\caption{\small Graphs based on numerical approximations of eigenfunctions corresponding to the first positive eigenvalue of the mixed  boundary value problem (above) and the Neumann one (below), for two domains $G_\varepsilon$ like those in Figures \ref{ej_unionprismas}-\ref{figura_boveda}.}
    \label{figura_mixed}
\end{center}
\end{figure}

\begin{figure}[ht]
\begin{center}
\scalebox{0.15}{\includegraphics{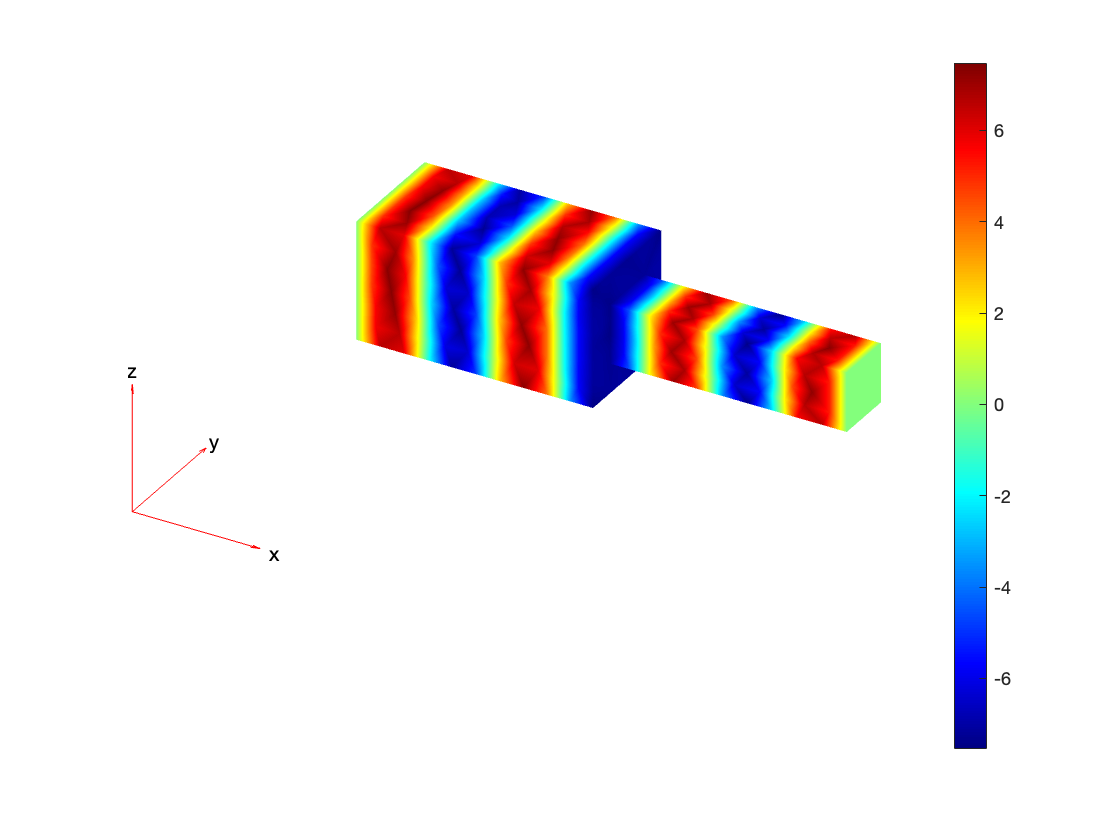}}
\scalebox{0.15}{\includegraphics{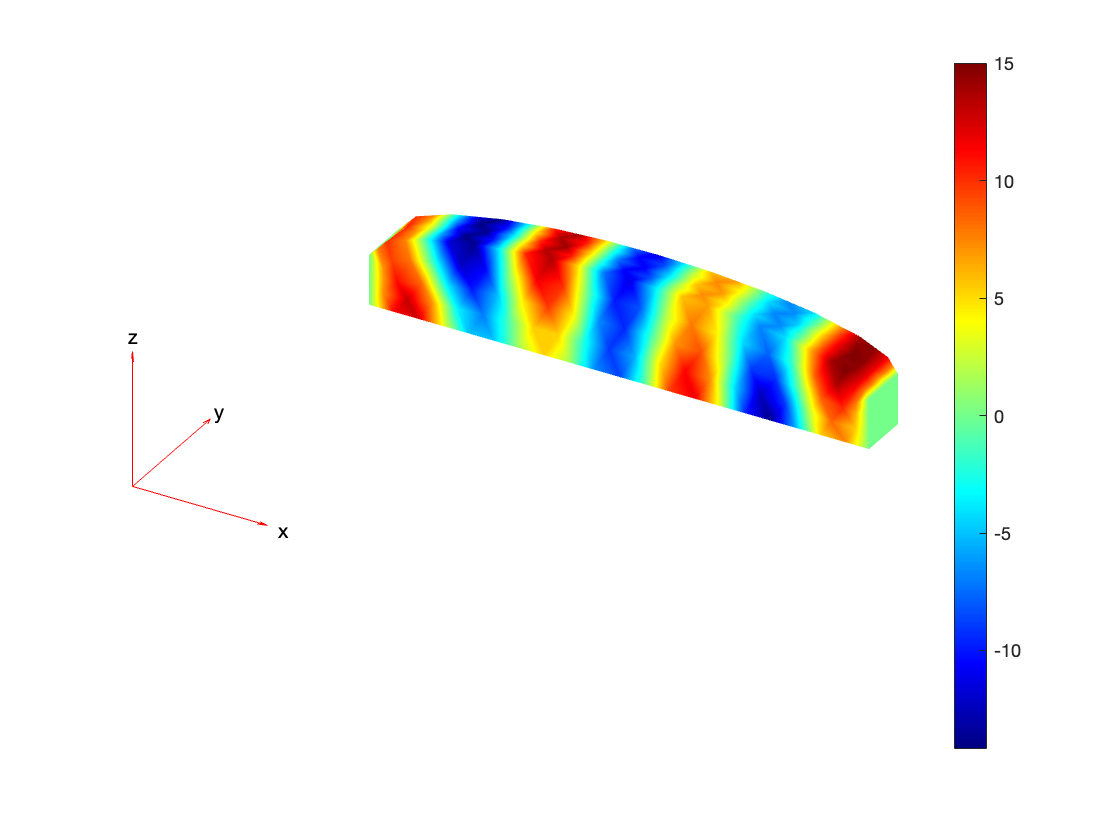}}
\scalebox{0.15}{\includegraphics{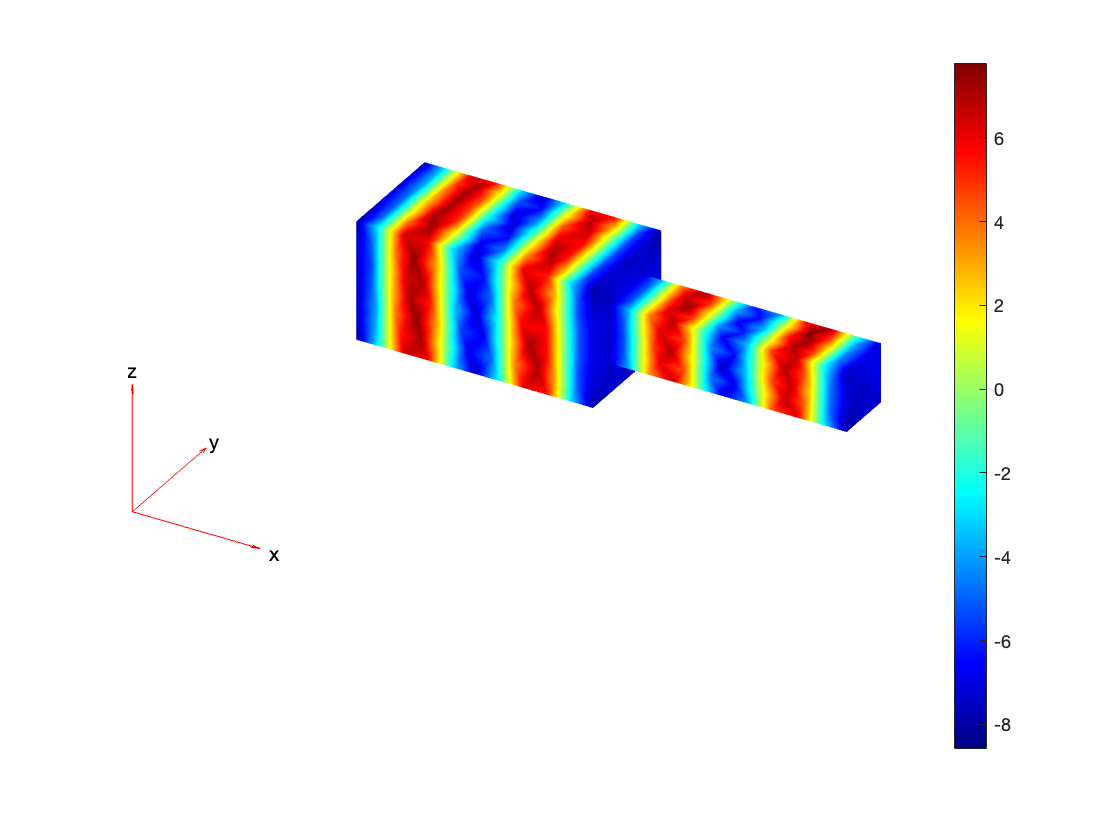}}
\scalebox{0.15}{\includegraphics{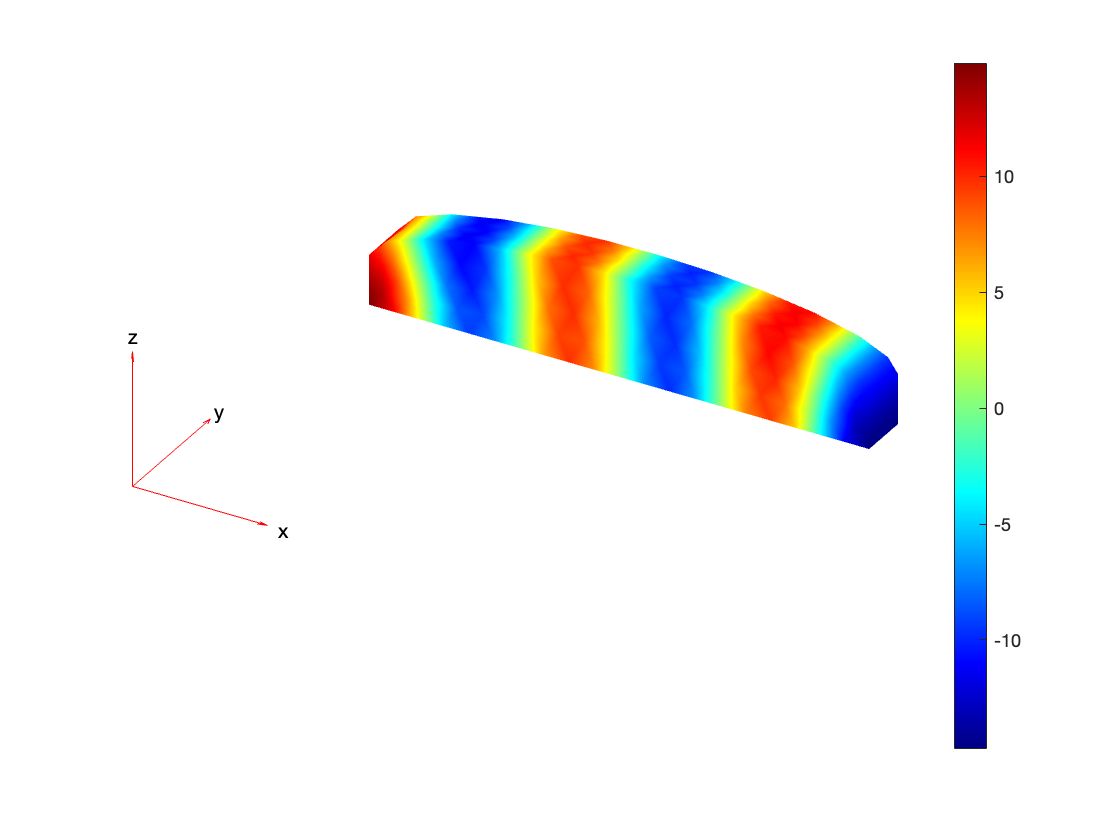}}
\caption{\small Graphs based on numerical  approximations of some eigenfunctions of the mixed boundary value problem (above) and the Neumann one (below), for two domains $G_\varepsilon$ like those in Figures \ref{ej_unionprismas}-\ref{figura_boveda}.}
    \label{figura_mixed2}
\end{center}
\end{figure}

\medskip
\noindent
To make it  easier for the reader, we provide Figures \ref{ej_unionprismas}-\ref{figura_boveda} as illustrations of samples of domains $G_\ee$ under consideration, and Figures \ref{figura_mixed}-\ref{figura_mixed2} as examples of eigenfunctions which present longitudinal oscillations, see \cite{articulo} for further precision on these geometrical configurations of the domains as well as for details and more graphics of  eigenfunctions.


\noindent
It should be emphasized that in Figure \ref{ej_unionprismas} many other polyhedral domains, union of rectangular prims or cuboids,   which satisfy the   required conditions for  $G_\ee$ and $G$, cf. \eqref{dominioG}-\eqref{dominioGeps},   could  be considered, always avoiding  cracks or cuspidal points,  cf. \cite{Mazya, Disser, Nazarov_book}. Also, the   dome-like top or wavy-like top domains in Figure \ref{figura_boveda} are  particular cases of more general domains whose expression can be
$$G_\varepsilon=\{(x_1,x_2,x_2) : x_1\in(l_0,l_1),
x_2\in(-\varepsilon h_1(x_1),\varepsilon h_2(x_1)),
x_3\in (-\varepsilon h_3(x_1),\varepsilon h_4(x_1))\},$$
where the functions $h_j$ are uniformly bounded regular functions subjects to conditions  \eqref{dominioG}-\eqref{acotacionarea}.
The same  can be said for the truncated cones or other different revolution domains without singular points. For particular varying cross sections and the elasticity system we refer to \cite{Viagno, Griso, NaSlTa} and references therein.

\medskip
\noindent
 Since the technique to show the  convergence of the eigenelements $(\lambda_\ee, u_\ee)$ involves a re-scaling of the domain $G_\varepsilon$  into a unit domain $G$, and henceforth introducing stretching variables,  in Section \ref{sec:ben2.1} we introduce some preliminary results that will be used through the rest of the paper.

\section{Preliminary results}\label{sec:ben2.1}

In this section, we introduce the variational formulation of problems  \eqref{formulacionuno_a},\eqref{formulacionuno_b},\eqref{fD} and  \eqref{formulacionuno_a},\eqref{formulacionuno_b},\eqref{fN} in the suitable Sobolev spaces as well as the reformulation of the problems in the fixed domain $G$ obtained by a change of scales in the transverse variables $x_2,x_3$, cf. \eqref{cambioGe-G}.

\medskip
\noindent
 The weak formulation of \eqref{formulacionuno_a},\eqref{formulacionuno_b},\eqref{fD}   reads: find $(\lambda_\ee, u_\ee)\in \mathbb{R}  \times H^1 (G_\ee,\Gamma_\ee^D)$, $u_\ee\not\equiv 0$,  satisfying
\begin{equation}\label{formulaciondebilproblemamixto}
 \int_{G_\ee} \nabla u_\ee.\nabla v  \, dx=\lambda_\ee  \int_{G_\ee} u^\ee v \,dx,  \quad\forall\, v\in H^1 (G_\ee,\Gamma_\ee^D),
\end{equation}
where $H^1 (G_\ee,\Gamma_\ee^D)$ denotes the  space completion of
\begin{equation*}
    \{u\in\mathscr{C}^\infty(\overline{G}_\varepsilon): u=0\text{ on }\Gamma_\varepsilon^D\},
\end{equation*}
equipped      with the norm generated by the scalar product $$(\nabla u,\nabla v)_{L^2(G_\ee)}.$$
On account of the Poincar\'e inequality, this norm is equivalent to the usual one in $H^1(G_\ee)$.

\noindent
The  formulation \eqref{formulaciondebilproblemamixto}  is classical in the couple of Hilbert spaces $H^1 (\Omega,\Gamma_\ee^D)\subset L^2(G_\ee)$ with a dense and compact embedding (cf. e.g. Section I.5 of \cite{SHSP89}), and therefore  the problem has a discrete spectrum.

\medskip
\noindent
For each fixed $\ee>0$,  let us denote   by
$$ 0<\lambda^1_\ee \le\lambda_\ee^2 \le\cdots \lambda_\ee^n \le \dots\to \infty, \quad \mbox{ as } n\to \infty , $$ the increasing sequence of eigenvalues, where we have adopted the convention of repeated eigenvalues according to their multiplicities.
In fact, using the minimax principle, we have shown in \cite{articulo} that  for each $n\in \mathbb{N}$ and for sufficiently small $\varepsilon>0$,  the eigenvalues satisfy the uniform bound
\begin{equation}\label{eq:ben8}
0 < C\leq \lambda_\varepsilon^n\leq C_n,
\end{equation}
where    $ C $ and $C_n$  are constants independent of $\ee$.
Also, we consider the corresponding set  of eigenfunctions $\{u_\ee^n\}_{ n=1}^\infty$ that  can be chosen to form an orthogonal basis in  $H^1 (G_\ee,\Gamma_\ee^D)$ and     in $L^2(G_\ee)$, subject to the normalization condition
\begin{equation}\label{eq:ben6}
\int_{G_\ee} \vert u_\ee\vert^2\, dx =\ee^2.
\end{equation}

\medskip
\noindent
Let us introduce  an auxiliary variable $y$, the so-called {\em   stretching variable}.
Its connection with $x$ is given by a change of variable  which transforms $G_\ee$ into $G$, namely,
\begin{equation}\label{cambioGe-G}
y_1= x_1, \quad y_2=\frac{x_2}{\ee}, \quad y_3=\frac{x_3}{\ee}.
\end{equation}

\noindent
Denoting $U_\varepsilon(y):=u_\varepsilon(y_1,\varepsilon y_2,\varepsilon y_3)$, we rewrite the problem \eqref{formulaciondebilproblemamixto} as
\begin{equation}
\int_G(\partial_{y_1}U_\varepsilon\partial_{y_1}V+\varepsilon^{-2}\partial_{y_2}U_\varepsilon\partial_{y_2}V+\varepsilon^{-2}\partial_{y_3}U_\varepsilon \partial_{y_3}V)dy=\lambda_\varepsilon\int_G U_\varepsilon Vdy\quad \forall V\in H^1(G,\Gamma^D),
    \label{formulaciondebilproblemamixtoescalado}
\end{equation}
and the normalization \eqref{eq:ben6} reads
\begin{equation}\label{eq:ben6b}
 \int_{G } \vert U_\ee\vert^2\, dy =1.
\end{equation}
Above, $\Gamma^D=\Gamma_0\cup\Gamma_1$ and  $H^1(G,\Gamma^D)$ stands for the completion of
\begin{equation*}
    \{U\in\mathscr{C}^\infty(\overline{G}): U=0\text{ on }\Gamma^D\}
\end{equation*}
with the  gradient norm which, by the Poincar\'e inequality,  is equivalent to  the $H^1$ norm in $G$.

\medskip
\noindent
Taking limits in \eqref{formulaciondebilproblemamixtoescalado} for particular test functions independent on the transverse variables, in  \cite{articulo}, we  have identified the limit problem of \eqref{formulaciondebilproblemamixto} with
\begin{equation}
    \int_{l_0}^{l_1} |D_{y_1}| \, \partial_{y_1}U_0 \varphi' \, dy_1 = \lambda_0 \int_{l_0}^{l_1} |D_{y_1}| \, U_0 \varphi \, dy_1 \quad \forall \varphi \in H_0^1(l_0,l_1),
    \label{formulaciondebilproblemalimite}
\end{equation}
which corresponds to the weak formulation of the following 1D Dirichlet problem \eqref{limitEDP},\eqref{limitD} and has a discrete spectrum.
Let
$$ 0<\lambda_0^1 \le\lambda_0^2 \le\cdots \lambda_0^n \le \dots\to \infty, \quad \mbox{ as } n\to \infty , $$ the increasing sequence of eigenvalues of \eqref{formulaciondebilproblemalimite} with the convention of repeated eigenvalues according to their multiplicities.

\noindent
As a matter of fact, the convergence result in \cite{articulo}  can be stated as follows:
If
\begin{equation}\label{benc}\lambda_\varepsilon  \to \lambda_0^*, \quad U_\varepsilon  \to U_0^* \quad \text{weakly in } H^1(G) \text{ as } \varepsilon \to 0,
\end{equation}
    then, $U_0^*$ depends only on the $y_1$ variable,  $\lambda_0^*$ is an eigenvalue of \eqref{formulaciondebilproblemalimite} and $U_0^*$ is an associated eigenfunction.
Let us observe that on account of \eqref{eq:ben6b}, the normalization \eqref{eq:ben6} along with estimates \eqref{eq:ben8} allow us to extract converging subsequences satisfying \eqref{benc}.

\medskip
\noindent
Let us proceed in a similar way with the Neumann problem \eqref{formulacionuno_a},\eqref{formulacionuno_b},\eqref{fN} which reads
\begin{equation}
    \left\{
    \begin{array}{ll}
        -\Delta u_\varepsilon=\lambda_\varepsilon u_\varepsilon & \text{ in }G_\varepsilon, \\[0.4em]
        \dfrac{\partial u_\varepsilon}{\partial n}=0 & \text{ on } \partial G_\varepsilon.
    \end{array}
    \right.
    \label{formulacionproblemaNeumann}
\end{equation}
\noindent
The weak formulation of problem \eqref{formulacionproblemaNeumann} reads: find  $(\lambda_\varepsilon,u_\varepsilon)\in\mathbb{R}\times H^1(G_\varepsilon)$ with $u_\varepsilon\not\equiv 0$ such that
\begin{equation}
    \int_{G_\varepsilon} \nabla u_\varepsilon\cdot\nabla vdx=\lambda_\varepsilon\int_{G_\varepsilon} u_\varepsilon vdx\quad \forall v\in H^1(G_\varepsilon).
    \label{formulaciondebilproblemaNeumann}
\end{equation}
As is well known, the first eigenvalue is $\lambda_\ee=0$ with associated eigenfunctions the constants. Therefore, we perform a shift to the spectrum  $\tilde{\lambda}_\varepsilon=\lambda_\varepsilon+1$ and
instead of \eqref{formulaciondebilproblemaNeumann} we consider the spectral problem
\begin{equation}
    \int_{G_\varepsilon} \nabla u_\varepsilon\cdot\nabla vdx+\int_{G_\varepsilon} u_\varepsilon vdx=\tilde{\lambda}_\varepsilon\int_{G_\varepsilon} u_\varepsilon vdx\quad \forall v\in H^1(G_\varepsilon),
    \label{formulaciondebilproblemaNeumann2}
\end{equation}
which has a discrete spectrum.   Its eigenvalues  satisfy
$$1=\tilde{\lambda}_\varepsilon^1<\tilde{\lambda}_\varepsilon^2\leq \dots\leq \tilde{\lambda}_\varepsilon^n\leq \dots\to\infty\quad \text{ when }n\to\infty$$
where we have assumed that they repeat according to their multiplicities.
Moreover, we  consider the set of eigenfunctions $\{u_\varepsilon^n\}_{n=1}^\infty$ to form an orthogonal basis of $H^1(G_\varepsilon)$ and $L^2(G_\varepsilon)$, subject to the following normalization condition:
\begin{equation}
    \int_{G_\varepsilon}|u_\varepsilon|^2dx=\varepsilon^2.
    \label{normalizacionproblemaNeumann}
\end{equation}
Using \eqref{cambioGe-G} and denoting $U_\varepsilon(y):=u_\varepsilon(y_1,\varepsilon y_2,\varepsilon y_3)$, we rewrite the problem \eqref{formulaciondebilproblemaNeumann2} as
\begin{equation}
\begin{split}
&\int_G(\partial_{y_1}U_\varepsilon\partial_{y_1}V+\varepsilon^{-2}\partial_{y_2}U_\varepsilon\partial_{y_2}V+\varepsilon^{-2}\partial_{y_3}U_\varepsilon \partial_{y_3}V)dy+\int_G U_\varepsilon Vdy \\
&\quad = \tilde{\lambda}_\varepsilon \int_G U_\varepsilon Vdy\quad\forall V\in H^1(G)
\end{split}
    \label{formulaciondebilproblemamixtoescaladoNeumann}
\end{equation}
and the normalization \eqref{normalizacionproblemaNeumann} becomes
\begin{equation*}
    \int_G |U_\varepsilon|^2dy=1.
\end{equation*}
Using a procedure similar to that of \cite{articulo}, we obtain the spectral limit problem of \eqref{formulaciondebilproblemaNeumann2}
\begin{equation}
    \int_{l_0}^{l_1} |D_{y_1}| \, \partial_{y_1}U_0 \varphi'  dy_1 +\int_{l_0}^{l_1}|D_{y_1}|U_0\varphi dy_1 = \tilde{\lambda}_0 \int_{l_0}^{l_1} |D_{y_1}| \, U_0 \varphi \, dy_1 \quad \forall \varphi \in H^1(l_0,l_1),
    \label{formulaciondebilproblemalimiteNeumann}
\end{equation}
which  also has a discrete spectrum. $\tilde{\lambda}_0 $  stands for the spectral parameter  and $U_0\in H^1(G)$ the corresponding eigenfunction. \eqref{formulaciondebilproblemalimiteNeumann}  corresponds to the weak formulation of the Neumann problem:
\begin{equation}
    \left\{
    \begin{array}{l}
         -\partial_{x_1}(|D_{x_1}| \partial_{x_1} U_0) + |D_{x_1}| U_0 = \tilde{\lambda}_0 |D_{x_1}| U_0, \quad x_1 \in (l_0,l_1), \\[0.3em]
         U_0'(l_0) = 0, \quad U_0'(l_1) = 0.
    \end{array}
    \right.
    \label{formulacionclasicaproblemalimiteNeumann}
\end{equation}
That is, \eqref{formulaciondebilproblemalimiteNeumann} is the weak formulation of \eqref{limitEDP},\eqref{limitN}, once we have introduced the shift of the spectral parameter $\tilde{\lambda}_0 = \lambda_0 +1$.
Let
$$ 1=\tilde{\lambda}_0^1 <\tilde{\lambda}_0^2 \le\cdots \tilde{\lambda}_0^n \le \dots\to \infty, \quad \mbox{ as } n\to \infty , $$ the increasing sequence of eigenvalues of \eqref{formulaciondebilproblemalimiteNeumann} with the convention of repeated eigenvalues according to their multiplicities.

\medskip
\noindent
Finally, for the sake of completeness let us introduce the result from the spectral perturbation theory, that will be used to prove the convergence of the above spectra with conservation of the multiplicity;  see Theorems III.1.4 and III.1.7 in  \cite{OlShYo92} for the proof.

 \begin{lemma}\label{LemaOl}
    Let $\mathcal{H}_\eps$ and $\mathcal{H}_0$ be two separable Hilbert spaces with the scalar products $\langle\cdot, \cdot\rangle_{\mathcal{H}_\eps}$ and $\langle\cdot, \cdot\rangle_{\mathcal{H}_0}$ respectively. Let $\mathcal{A}_\eps\in \mathcal{L}(\mathcal{H}_\eps)$ and  $\mathcal{A}_0\in \mathcal{L}(\mathcal{H}_0)$. Let $\mathcal{V}$ be a subspace of $\mathcal{H}_0$ such that $Im(A_0)= \{v: \, v= A_0 u \: \mbox{ \rm with } \: u\in \mathcal{H}_0\}\subset \mathcal{V}$. We assume that the following properties are satisfied:
    \begin{enumerate}
    \item[$C_1)$] There exist continuous linear operators $R_\varepsilon:\mathcal{H}_0\rightarrow\mathcal{H}_\varepsilon$ and a constant $\gamma>0$ such that
    \begin{equation}
        \langle R_\varepsilon f_0, R_\varepsilon f_0 \rangle_{\mathcal{H}_\varepsilon}\to \gamma\langle f_0, f_0 \rangle_{\mathcal{H}_0}\quad\text{ when }\varepsilon\to 0\quad\text{ for any }f_0\in\mathcal{V}.
        \label{C1}
    \end{equation}

    \item[$C_2)$] The operators $\mathcal{A}_\varepsilon$ and $\mathcal{A}_0$ are positive, compact, and self-adjoint. Moreover, the norms of the operators $\mathcal{A}_\varepsilon$ are bounded by a constant independent of $\varepsilon$.

    \item[$C_3)$] For any $f\in\mathcal{V}\subset\mathcal{H}_0$ it holds that
    \begin{equation}
        ||\mathcal{A}_\varepsilon R_\varepsilon f-R_\varepsilon \mathcal{A}_0 f||_{\mathcal{H}_\varepsilon}\to 0\quad\text{ when }\varepsilon\to 0.
       \label{C3}
    \end{equation}

    \item[$C_4)$] The family of operators $\{\mathcal{A}_\varepsilon\}_{\varepsilon>0}$ is uniformly compact in the following sense: for each sequence $\{f_\varepsilon\}_{\varepsilon>0}\text{ in }\mathcal{H}_\varepsilon$ such that $\sup_{\varepsilon}||f_\varepsilon||_{\mathcal{H}_\varepsilon}<\infty$, we can extract a subsequence $\{f_{\varepsilon'}\}_{\varepsilon'>0}$ such that for some $w_0\in\mathcal{V}$ it holds that
    \begin{equation*}
        ||\mathcal{A}_{\varepsilon'}f_{\varepsilon'}-R_{\varepsilon'}w_0||_{\mathcal{H}_{\varepsilon'}}\to 0\quad\text{ when }\varepsilon'\to 0.
    \end{equation*}

\end{enumerate}
Let $\{\mu_\eps^i\}_{i=1}^\infty$ ($\{\mu_0^i\}_{i=1}^\infty$, respectively) be the sequence of the eigenvalues of $\mathcal{A}_\eps$ ($\mathcal{A}_0$, respectively) with the usual convention of repeated eigenvalues. Let $\{u_\eps^i\}_{i=1}^\infty$ ($\{u_0^i\}_{i=1}^\infty$, respectively) be the corresponding eigenfunctions which is assumed
to be an orthonormal basis in $\mathcal{H}_\varepsilon$ ($\mathcal{H}_0$, respectively).
\\
Then, for each fixed $k$, we have
$$\mu_\eps^k \to \mu_0^k\quad \mbox{ as } \quad \eps\to 0\,.$$
Moreover, for each sequence we can extract a subsequence $\eps'\to 0$ such that for any $k$,
$$
\|A_{\eps '} u_{\eps'}^k- {\mathcal R}_{\eps'} u_*^k\|_{\mathcal{H}_{\eps'}} \to 0\quad \mbox{ as } \quad \eps'\to 0\,,$$
where $u_*^k$ is an eigenfunction of $\mathcal{A}^0$ associated with $\mu_0^k$ and the set $\{u_*^i\}_{i=1}^\infty$ forms an orthogonal basis of $\mathcal{H}_0$.
\\
In addition, for each  eigenfunction $u_0$ associated to $\mu_0^{k+1}$ eigenvalue with multiplicity $m$, $\mu_0^{k}>\mu_0^{k+1}=\mu_0^{k+2}=\dots=\mu_0^{k+m}>\mu_0^{k+m+1}$,  with $||u_0||_{\mathcal{H}_0}=1$, there exists a linear combination $\bar{u}_\varepsilon$ of the eigenfunctions $u_\varepsilon^{k+1},\dots,u_\varepsilon^{k+m}$ such that
    \begin{equation*}
        ||\bar{u}_\varepsilon-u_0||_{\mathcal{H}_\varepsilon}\to 0\quad\text{ when }\varepsilon\to 0.
    \end{equation*}
\label{lemaole}
\end{lemma}
\noindent It could be noted that the last result in the lemma is proved by Theorem III.1.7 in  \cite{OlShYo92}. In the case where $m=1$, $\bar{u}_\varepsilon$ is an  eigenfunction corresponding to  $\mu_\varepsilon^{k+1}$.

\section{\!\!The convergence for the mixed boundary value problem}\label{sec:ben3}

In this section,  we prove convergence of the eigenvalues \eqref{formulacionuno_a},\eqref{formulacionuno_b},\eqref{fD}   towards those of   \eqref{limitEDP},\eqref{limitD}  with conservation of the multiplicity (cf, \eqref{conservation}). To do this, we consider the problems in the  stretching variable \eqref{cambioGe-G} and    use Lemma \ref{LemaOl}. Hence, we are led to show the convergence of solutions, as $\ee\to 0$,  for stationary problems associated to \eqref{formulaciondebilproblemamixtoescalado},   towards the solutions  of stationary problems associated to \eqref{formulaciondebilproblemalimite}.

\noindent
For problem   \eqref{formulaciondebilproblemamixtoescalado},   let us introduce   the  spaces and operators under consideration $\mathcal{H}_\varepsilon$, $\mathcal{H}_0$, $\mathcal{A}_\varepsilon$, $\mathcal{A}_0$ and $R_\varepsilon$ and verify conditions $C_1-C_4$ arising in the statement of Lemma \ref{lemaole}.
\begin{itemize}
    \item {\bf Spaces $\mathcal{H}_\varepsilon$ and $\mathcal{V}_\varepsilon$:}   Fixed $\ee$, we consider $\mathcal{H}_\varepsilon=L^2(G)$   and the space $\mathcal{V}_\varepsilon$ to be  $ H^1(G,\Gamma^D)$ endowed with the scalar product given by the bilinear form
        \begin{equation}\label{formaPhie}
       \langle U,V\rangle_{\mathcal{V}_\ee} :=  \int_G(\partial_{y_1}U\partial_{y_1}V+
       \varepsilon^{-2}\partial_{y_2}U\partial_{y_2}V
       +\varepsilon^{-2}\partial_{y_3}U\partial_{y_3}V)dy.
       \end{equation}
        Considering the Poincar\'e inequality, the norm generated by \eqref{formaPhie} is equivalent to the usual one in $H^1(G)$.

    \item {\bf Spaces $\mathcal{H}_0$ and $\mathcal{V}_0$:}  Let $\mathcal{H}_0$ be the space $L^2(l_0,l_1)$ with the weighted  scalar product
 \begin{equation}\label{w1}
\langle u,v\rangle_{\mathcal{H}_0}:=\int_{l_0}^{l_1} |D_{x_1}| uvdx_1.
 \end{equation}
    Let   $\mathcal{V}_0$ be    the space $H_0^1(l_0,l_1)$ endowed with the scalar product given by
   \begin{equation}\label{formaPhi0}\langle u,v\rangle_{\mathcal{V}_0}:=\int_{l_0}^{l_1} |D_{x_1}| u'v'dx_1. \end{equation}
    Also, the norm associated with \eqref{formaPhi0}  in $\mathcal{V}_0$ is equivalent to the usual one in $H^1(l_0,l_1)$ due to the Poincar\'e inequality  and the bound \eqref{acotacionarea}.

    \item {\bf Operator $\mathcal{A}_\varepsilon$:} We define  operator $\mathcal{A}_\varepsilon:\mathcal{H}_\varepsilon\to\mathcal{H}_\varepsilon$ as:  $\mathcal{A}_\varepsilon f=u_f^\varepsilon$, where $f\in\mathcal{H}_\varepsilon$ and $u_f^\varepsilon\in\mathcal{V}_\varepsilon\subset\mathcal{H}_\varepsilon$ is the unique solution of the problem
    \begin{equation}\Phi_\varepsilon(u_f^\varepsilon,v)=\langle f, v\rangle_{\mathcal{H}_\varepsilon}\quad\forall v\in\mathcal{V}_\varepsilon,  \label{ProblemaOperadorAe}
    \end{equation}
    $\Phi_\ee$ being the form defined by the scalar product \eqref{formaPhie}, namely,
    \begin{equation}\label{formaPhiefi}
       \Phi_\ee(U,V) :=\langle U,V\rangle_{\mathcal{V}_\ee} \quad \forall U,V\in  {\mathcal{V}_\ee}. \end{equation}
        The existence and uniqueness of the solution  of   \eqref{ProblemaOperadorAe} follows  from the Lax-Milgram Lemma (cf. \eqref{formaPhie} and \eqref{formaPhiefi}).

     \noindent
     Next, let us show that $\mathcal{A}_\varepsilon$ is a continuous, compact, self-adjoint, and positive operator.

   \noindent
   {\em Continuity:} Using the Poincaré inequality in $H^1(G,\Gamma^D)$,    we have
   $$||u_f^\varepsilon||_{L^2(G)}\leq C||\nabla u_f^\varepsilon||_{L^2(G)}.$$
        This inequality along with  \eqref{ProblemaOperadorAe} (cf. \eqref{formaPhie} and \eqref{formaPhiefi}), and the Cauchy-Schwarz inequality, gives the chain of inequalities
        $$||\nabla u_f^\varepsilon||_{L^2(G)}^2\leq \Phi_\varepsilon(u_f^\varepsilon,u_f^\varepsilon)=\langle f,u_f^\varepsilon\rangle_{\mathcal{H}_\varepsilon}
         {\leq} ||f||_{\mathcal{H}_\varepsilon}||u_f^\varepsilon||_{\mathcal{H}_\varepsilon} {\leq} C||f||_{\mathcal{H}_\varepsilon}||\nabla u_f^\varepsilon||_{\mathcal{H}_\varepsilon}.$$

        \noindent
        Hence,
        \begin{equation}
            ||\nabla u_f^\varepsilon||_{L^2(G)}\leq C||f||_{L^2(G)},
            \label{cotagradiente}
        \end{equation}
        and therefore, using again the Poincaré inequality, we have
         $||u_f^\varepsilon||_{L^2(G)}\leq C||f||_{L^2(G)},  $
         with $C$ a constant independent of $f$ and $\varepsilon$.
         In this way, we have obtained the uniform bound
       \begin{equation}\label{unif} ||\mathcal{A}_\varepsilon f||_{L^2(G)}\leq C||f||_{L^2(G)},\end{equation}
     which implies  the continuity of the operator $\mathcal{A}_\varepsilon$ in  ${\mathcal{H}_\varepsilon} $.   Moreover, it  implies that   the set $\{  \Vert \mathcal{A}_\varepsilon\Vert_{\mathcal L(\mathcal{H}_\varepsilon)} \}_{\varepsilon}$    is uniformly bounded.

  \noindent
  {\em Compactness:}
       Let us show that if $f_n\to f$ in $L^2(G)$-weak as $n\to\infty$, then $u_{f_n}^\varepsilon=\mathcal{A}_\varepsilon f_n\to u_f^\varepsilon=\mathcal{A}_\varepsilon f$   in $L^2(G)$ as $n\to\infty$.
       By the definition of $\mathcal{A}_\varepsilon$, we have that $u_{f_n}^\varepsilon, u_f^\varepsilon\in\mathcal{V}_\varepsilon$, and   using the bound  \eqref{cotagradiente}  for $f \equiv f_n$ , we establish the bound of
       \begin{equation}\label{c1}
       ||\nabla u_{f_n}^\varepsilon||_{L^2(G)}\leq C
       \end{equation}
       where $C$ is a constant independent of $n$. Hence, for any sequence of $n$, we can extract a subsequence, still denoted by $n$, such that $u_{f_n}^\varepsilon$ converges weakly in    $H^1(G)$    towards some function $u_*^\varepsilon$,  as $n\to\infty$. Once that we identify  $u_*^\varepsilon$ with $u_f^\varepsilon$ the solution of \eqref{ProblemaOperadorAe}, by the uniqueness of solution, we deduce the weak convergence of the whole sequence in $H^1(G)$ and,  since the embedding  $\mathcal{V}_\varepsilon=H^1(G,\Gamma^D)$ into $\mathcal{H}_\varepsilon=L^2(G)$ is compact, the above convergence for  $u_{f_n}^\varepsilon$  towards $u_{f}^\varepsilon$ holds in  the strong topology of $L^2(G)$.

        \noindent
        To show the  equality $u_*^\varepsilon=u_f^\varepsilon$ it suffices to take limits in  $$
        \int_G (\partial_{y_1}u_{f_n}^\varepsilon\partial_{y_1}v+\varepsilon^{-2}\partial_{y_2}u_{f_n}^\varepsilon\partial_{y_2}v+\varepsilon^{-2}\partial_{y_3}u_{f_n}^\varepsilon \partial_{y_3}v)dy=\int_G f_nvdy\quad\forall v\in\mathcal{V}_\varepsilon.
     $$
      as $n\to \infty$ and get equation \eqref{ProblemaOperadorAe}.

\noindent{\em Self-adjointness:} Since the bilinear form  $\Phi_\ee$ is symmetric, we can write
\begin{equation*}
\begin{split}
&\langle \mathcal{A}_\varepsilon f,g\rangle_{\mathcal{H}_\varepsilon}
=\langle u_f^\varepsilon,g\rangle_{\mathcal{H}_\varepsilon}
=\langle g,u_f^\varepsilon\rangle_{\mathcal{H}_\varepsilon}
=\Phi_\varepsilon(u_g^\varepsilon,u_f^\varepsilon)\\
&\quad
=\Phi_\varepsilon(u_f^\varepsilon,u_g^\varepsilon)
=\langle f,u_g^\varepsilon\rangle_{\mathcal{H}_\varepsilon}
=\langle f,\mathcal{A}_\varepsilon g\rangle_{\mathcal{H}_\varepsilon}.
\end{split}
\end{equation*}

\noindent
{\em Positiveness:} This is straightforward because $\langle \mathcal{A}_\varepsilon f,f\rangle=\Phi_\varepsilon(u_f^\varepsilon,u_f^\varepsilon)\geq 0$ and it only vanishes if $u_f^\varepsilon\equiv0$, cf. \eqref{formaPhiefi},  which implies $f\equiv0$ by the uniqueness of the solution of problem \eqref{ProblemaOperadorAe}.

\item {\bf Operator $\mathcal{A}_0$:} The operator $\mathcal{A}_0:\mathcal{H}_0\to\mathcal{H}_0$  defined as  $\mathcal{A}_0 f=u_f$, for  $f\in\mathcal{H}_0$  and  $u_f\in\mathcal{V}_0\subset\mathcal{H}_0$  the unique solution of the problem
        \begin{equation*}
        \Phi_0(u_f,v)=\langle f, v\rangle_{\mathcal{H}_0} \quad \forall v \in \mathcal{V}_0,
    \end{equation*}
    where $\Phi_0 $ is defined by the scalar product in \eqref{formaPhi0};
    namely,
$$\Phi_0(u,v):=  \langle u, v\rangle_{\mathcal{V}_0} \quad \forall u , v \in {\mathcal{V}_0}.$$
 Note that the existence and uniqueness of the solution $u_f$ follows again from Lax-Milgram   Lemma.
    \\
   The  fact that the operator $\mathcal{A}_0$ is continuous, compact, self-adjoint   and positive follows  using arguments similar to those for the operator $\mathcal{A}_\varepsilon$.

    \item {\bf Operator $R_\varepsilon$:} Let us define the operator $R_\varepsilon:\mathcal{H}_0\to\mathcal{H}_\varepsilon$  as an extension operator independent of $\varepsilon$. Namely, for $f\in \mathcal{H}_0=L^2(l_0,l_1)$,
        \begin{equation}\label{ext}
         (R_\varepsilon f)(x_1,x_2,x_3):=f(x_1).
         \end{equation}
 \end{itemize}

\medskip
\noindent
Now, let us   verify conditions $C_1-C_4$ in Lemma \ref{lemaole}.
\begin{itemize}
    \item {\bf Verification of $C_1$:}  
   For $f_0\in\mathcal{V}$, writing
  \begin{equation*}
  \begin{split}
  &||R_\varepsilon f_0||_{\mathcal{H}_\varepsilon}^2
  =\int_{G}|f_0|^2dy
  =\int_{l_0}^{l_1}|f_0(x_1)|^2dx_1\int_{D_{x_1}}dy_2dy_3\\
  &\quad =\int_{l_0}^{l_1}|D_{x_1}|\text{ }|f_0(x_1)|^2dx_1
  =||f_0||_{\mathcal{H}_0}^2,
  \end{split}
  \end{equation*}
\eqref{C1} holds taking $\gamma=1^.$
    \item {\bf Verification of $C_2$:} It was already verified when introducing the operators $\mathcal{A}_\varepsilon$ and $\mathcal{A}_0$ that they are positive, compact, and self-adjoint and, moreover,  from \eqref{unif}
      \begin{equation}\label{unif1}
      \Vert \mathcal{A}_\varepsilon\Vert_{\mathcal L(\mathcal{H}_\varepsilon)} \leq  C,\end{equation}
      for a constant $C$ independent of $\varepsilon$.

    \item {\bf Verification of $C_3$:} Let us show that for any $f\in\mathcal{V}_0\subset\mathcal{H}_0$ we have \eqref{C3}.
    Fixing $f\in H_0^1(l_0,l_1)$, we denote $u_f^\varepsilon=\mathcal{A}_\varepsilon R_\varepsilon f$.  
    The goal is to prove that $u_f^\varepsilon\to u_f^0$ in $\mathcal{H}_\varepsilon=L^2(G)$ when $\varepsilon\to 0$, where $u_f^0$ satisfies
    \begin{equation}
        \Phi_0(u_f^0,v)=\langle f,v\rangle_{\mathcal{H}_0}\quad\forall v\in H_0^1(l_0,l_1).
        \label{Phi0ComprobacionC3}
    \end{equation}
    and hence, taking into account that $R_\ee$ is the extension operator \eqref{ext}, the convergence   $$||\mathcal{A}_\varepsilon R_\varepsilon f-R_\varepsilon \mathcal{A}_0 f||_{\mathcal{H}_\varepsilon}\to 0\quad\text{ when }\varepsilon\to 0, $$
   also holds.

\noindent
   To show the above convergence,  first, we   show  that $u_f^\varepsilon$  is  uniformly bounded in $H^1(G)$. Starting from the weak formulation
    \begin{equation}
        \Phi_\varepsilon(u_f^\varepsilon,v)=\langle R_\varepsilon f,v\rangle_{\mathcal{H}_\varepsilon}\quad\forall v\in H^1(G,\Gamma^D),
        \label{PhieComprobacionC3}
    \end{equation}
    we take $v=u_f^\varepsilon$ and  write
    $$\Phi_\varepsilon(u_f^\varepsilon,u_f^\varepsilon)=\langle f,u_f^\varepsilon\rangle_{\mathcal{H}_\varepsilon}.$$
    Recalling the definition \eqref{formaPhiefi} (cf. \eqref{formaPhie}), we have the chain of inequalities
    $$||\nabla u_f^\varepsilon||_{L^2(G)}^2\leq \Phi_\varepsilon(u_f^\varepsilon,u_f^\varepsilon)=\langle f,u_f^\varepsilon\rangle_{\mathcal{H}_\varepsilon} {\leq} ||f||_{L^2(G)}||u_f^\varepsilon||_{L^2(G)} {\leq} C||f||_{L^2(G)}||\nabla u_f^\varepsilon||_{L^2(G)},$$
    where  we have used the Cauchy-Schwarz and Poincar\'e inequalities and $C$ is a constant independent of $\varepsilon$.
    \\
   Consequently, we obtain
    \begin{equation}\label{nna} ||\nabla u_f^\varepsilon||_{L^2(G)}\leq C||f||_{L^2(G)}\quad  \text{ and } \quad ||  u_f^\varepsilon||_{L^2(G)}\leq C||f||_{L^2(G)}.\end{equation}
    This along with the trace   theorem ensures the existence of   converging subsequences of $\{u_f^\varepsilon\}_\varepsilon$  in  $H^1(G,\Gamma^D)$-weak.

    \noindent
    Also, from the above chain of inequalities  we deduce (cf. \eqref{formaPhie})
    \begin{equation}
         ||\partial_{y_2}u_f^\varepsilon||_{L^2(G)}\leq \hat{K} \varepsilon||f||_{L^2(G)} \text{ and } ||\partial_{y_3}u_f^\varepsilon||_{L^2(G)}\leq \hat{K} \varepsilon||f||_{L^2(G)}
        \label{convergenciaderivaday2}
    \end{equation}
    with $\hat K$ a constant independent of $\varepsilon$.
   Bounds in \eqref{convergenciaderivaday2}   imply  the convergence towards zero of the derivatives of $\partial_{y_2}u_f^\varepsilon$ and $\partial_{y_3}u_f^\varepsilon$ in $L^2(G)$ and a simple argument leads  us to assert  that  the limit functions of the subsequences $u_f^\varepsilon$ do not depend on the variables $y_2$ and $y_3$.
    \\
  From \eqref{nna},  $u_f^\varepsilon$ is uniformly bounded in $H^1(G)$  and, therefore, for each sequence there exists a subsequence $\varepsilon'$ such that $u_f^{\varepsilon'}\to u$ in $H^1(G)$-weak and  strong in $L^2(G)$ for some $u\in H^1(G,\Gamma^D)$.
    Next, let us   identify $u$ with $u_f^0=\mathcal{A}_0f$ and, then, by the uniqueness of the solution of problem \eqref{Phi0ComprobacionC3}, the whole sequence will converge to $u_f^0$.

    \noindent
    To this end, in \eqref{PhieComprobacionC3},  we take functions $v_\varepsilon\in H^1(G, \Gamma^D)$
     $$v_\varepsilon(y_1,y_2,y_3):=v(y_1)\quad\text{ with }v\in\mathscr{C}_0^\infty(l_0,l_1)$$ and write
    \begin{equation*}
        \int_G \partial_{y_1}u_f^\varepsilon\partial_{y_1}vdy = \langle R_\varepsilon f,v_\varepsilon\rangle_{L^2(G)}  .\end{equation*}
        Taking limits as $\ee \to 0$, since  by definition  
    $$\langle R_\varepsilon f,v_\varepsilon\rangle_{L^2(G)}=\langle f,v\rangle_{\mathcal{H}_0}, $$
    we get
      $$ \int_G \partial_{y_1}u\partial_{y_1}vdy=
         \int_{l_0}^{l_1}|D_{x_1}|\partial_{x_1}u\partial_{x_1}vdx_1=\Phi_0(u,v)=\langle f,v\rangle_{\mathcal{H}_0} .$$
  Thus, since $u=0$ on $\Gamma^D$, a density argument allows us to conclude that
    $$\Phi_0(u,v)=\langle f,v\rangle_{\mathcal{H}_0}\quad\forall v\in H_0^1(l_0,l_1),$$
    and therefore it satisfies \eqref{Phi0ComprobacionC3}. By the uniqueness of the limit, $u=\mathcal{A}_0f$. That is, $u_f^\varepsilon\to\mathcal{A}_0 f$  in $L^2(G)$, as $\varepsilon \to 0$, and on account of the definition of the extension operator, cf. \eqref{ext}, we have concluded the verification of condition $C_3$.

    \item {\bf Verification of $C_4$:}
    Let us show that for each sequence $\{f_\varepsilon\}_{\varepsilon>0}\text{ in }\mathcal{H}_\varepsilon$ such that $\sup_{\varepsilon}||f_\varepsilon||_{\mathcal{H}_\varepsilon}<\infty$, we can extract a subsequence $\{f_{\varepsilon'}\}_{\varepsilon'>0}$ such that for some $w_0\in\mathcal{V}_0$ we have
    $$||\mathcal{A}_{\varepsilon'}f_{\varepsilon'}-R_{\varepsilon'}w_0||_{\mathcal{H}_{\varepsilon'}}\to 0\quad\text{ when }\varepsilon'\to 0.$$

\noindent
    We consider the sequence $u_\varepsilon:=\mathcal{A}_\varepsilon f_\varepsilon\in H^1(G)$, i.e., $u_\varepsilon$ is the solution of
    \begin{equation*}
        \Phi_\varepsilon(u_\varepsilon,v)=\langle  f_\varepsilon,v\rangle_{\mathcal{H}_\varepsilon}\quad\forall v\in H^1(G,\Gamma^D).
    \end{equation*}
    Since the operators $\mathcal{A}_\varepsilon$ are uniformly bounded and the sequence $\{f_\varepsilon\}_{\varepsilon>0}$ is bounded, then the sequence $\{u_\varepsilon\}_{\varepsilon>0}$ defined above is bounded in $L^2(G)$, that is, $||u_\varepsilon||_{\mathcal{H}_\varepsilon}\leq C$ for some constant $C>0$. Then, we rewrite  the reasoning in  \eqref{PhieComprobacionC3}--\eqref{convergenciaderivaday2},  to  obtain the  boundeness   $|| u_\varepsilon||_{H^1(G)}\leq K$ and that the possible limits of $\{u_\varepsilon\}_{\varepsilon>0}$ are independent of $y_2$ and $y_3$.
    \\
    In this way, we can extract a subsequence $\{u_{\varepsilon'}\}_{\varepsilon'>0}$ converging strongly in $L^2(G)$ and weakly in $H^1(G)$ towards some $u\in H^1(G,\Gamma^D)$.
    Hence,  since $u=0$ on $\Gamma^D$ and it is independent of $y_2$ and $y_3$, we can identify $u$ with $R_{\varepsilon'} w_0$ for some $w_0\in\mathcal{V}_0=H_0^1(l_0,l_1)$, which concludes the verification of condition $C_4$.
\end{itemize}

\noindent
Finally, it is worth noting that the eigenvalues $\{\mu_\varepsilon^n\}_{n=1}^\infty$ of the operator $\mathcal{A}_\varepsilon$ are the inverses of the eigenvalues of problem \eqref{formulaciondebilproblemamixtoescalado}. Namely,
    \begin{equation}\label{i1}
        \mu_\varepsilon^n=\frac{1}{\lambda_\varepsilon^n}, \quad n=1,2,\cdots.
    \end{equation}
    Similarly, the eigenvalues $\{\mu_0^n\}_{n=1}^\infty$ of $\mathcal{A}_0 $  are the inverses of  the eigenvalues of \eqref{formulaciondebilproblemalimite}:
    \begin{equation}\label{i2}
        \mu_0^n=\frac{1}{\lambda_0^n}, \quad n=1,2,\cdots.
    \end{equation}

\medskip
\noindent
Gathering all the results above, we apply Lemma \ref{lemaole} and, considering \eqref{i1} and \eqref{i2},  we have proved the convergence of the spectrum of \eqref{formulacionuno_a},\eqref{formulacionuno_b},\eqref{fD} towards that of \eqref{limitEDP},\eqref{limitD} with conservation of the multiplicity as well as the convergence of the corresponding eigenfunctions  as we state in the following theorem.
\begin{theorem}\label{TheoremM}
Suppose that the  area of the  cross sections of the domain $G$, $|D_{x_1}|$,    is a piecewise continuous function in $ [l_0,l_1]$,  satisfying   \eqref{acotacionarea}. Let $\{\lambda_\varepsilon^n\}_{n=1}^\infty$ be the eigenvalues of problem \eqref{formulaciondebilproblemamixtoescalado} and $\{U_\varepsilon^n\}_{n=1}^\infty$ be the corresponding eigenfunctions forming an orthonormal basis in $L^2(G)$ (see \eqref{eq:ben6b}). Then, for any $n\in\mathbb{N}$, we have that
\begin{equation}\label{conservation}
\lambda_\varepsilon^n\to\lambda_0^n\quad\text{ when }\varepsilon\to 0,
\end{equation}
where $\lambda_0^n$ is the $n$-th eigenvalue of  \eqref{formulaciondebilproblemalimite}. Moreover, for each sequence $\varepsilon$, we can extract a subsequence (still denoted by $\varepsilon$) such that
\begin{equation}\label{conservationfp}U_\varepsilon^n\to U_*^n\quad\text{in }H^1(G)\text{-weak when }\varepsilon\to 0,\end{equation}
where $U_*^n$ is an eigenfunction of \eqref{formulaciondebilproblemalimite} associated with $\lambda_0^n$ and such that $\{U_*^n\}_{n=1}^\infty$ form a basis of $H_0^1(l_0,l_1)$ and $L^2(l_0,l_1)$.
\\
In addition, for each $U_0$ eigenfunction of the limit problem \eqref{formulaciondebilproblemalimite} associated with $\lambda_0^{k+1}$ eigenvalue with multiplicity $m$, $\lambda_0^{k}<\lambda_0^{k+1}=\lambda_0^{k+2}=\dots=\lambda_0^{k+m}<\lambda_0^{k+m+1}$, with $||U_0||_{\mathcal{H}_0}=1$, there exists a linear combination $\bar{U}_\varepsilon$ of the eigenfunctions $U_\varepsilon^{k+1},\dots,U_\varepsilon^{k+m}$ of problem \eqref{formulaciondebilproblemamixtoescalado} such that
    \begin{equation*}
        ||\bar{U}_\varepsilon-U_0||_{L^2(G)}\to 0\quad\text{ when }\varepsilon\to 0.
    \end{equation*}
\end{theorem}

\section{The convergence for the Neumann problem}\label{sec:ben4}
\noindent
In this section we address the spectral convergence for the Laplacian with Neumann boundary conditions, cf. \eqref{formulacionproblemaNeumann} or, equivalently, \eqref{formulaciondebilproblemaNeumann}. As outlined in Section  \ref{sec:ben2.1}, we  work with the shifted spectrum  $\tilde{\lambda}_\varepsilon=\lambda_\varepsilon+1$ and show the convergence of the eigenvalues of \eqref{formulaciondebilproblemamixtoescaladoNeumann}  towards those of  \eqref{formulaciondebilproblemalimiteNeumann} with conservation of the multiplicity (cf. \eqref{conservation2}).
We follow the technique  in Section \ref{sec:ben3}  with the appropriate modifications, that we introduce below. Among other things, due the above shift of the spectrum,  we do not need to use the Poincar\'e inequality and some steps can be simplified. For the sake of brevity we avoid repetitions as much as possible.

\medskip

\noindent
First, let us introduce the necessary spaces and operators $\mathcal{H}_\varepsilon$, $\mathcal{H}_0$, $\mathcal{A}_\varepsilon$, $\mathcal{A}_0$ and $R_\varepsilon$ to apply Lemma \ref{lemaole}.

\noindent
The spaces $\mathcal{H}_\varepsilon$ and $\mathcal{H}_0$  are those of Section \ref{sec:ben3}, namely $L^2(G)$ and $L^2(l_0,l_1)$ with  the norm generated by \eqref{w1}. Let  $\mathcal{V}_\varepsilon$ be the space $H^1(G)$   with the scalar product defined by:
$$\langle U,V\rangle_{\mathcal{V}_\ee} :=\int_G(\partial_{y_1}U\partial_{y_1}V+\varepsilon^{-2}\partial_{y_2}U\partial_{y_2}V+\varepsilon^{-2}\partial_{y_3}U\partial_{y_3}V)dy+\int_G U Vdy. $$
Let $\mathcal{V}_0$ be $H^1(l_0,l_1)$ with the scalar product defined by
$$\langle u,v\rangle_{\mathcal{V}_0}=\int_{l_0}^{l_1} |D_{x_1}| u'v'dx_1+\int_{l_0}^{l_1}|D_{x_1}|uvdx_1.$$

\noindent
On these spaces we introduce the forms $\Phi_\ee$ and $\Phi_0$ given by the scalar products, namely
$$\Phi_\ee(U ,V ):=   \langle U,V\rangle_{\mathcal{V}_\ee} \quad \text{ and } \quad \Phi_0(u ,v):= \langle u,v\rangle_{\mathcal{V}_0},$$
respectively. It is self-evident that they define bilinear, continuous, symmetric and coercive forms on ${\mathcal{V}_\ee}$ and ${\mathcal{V}_0}$ respectively.

\noindent
As for the operator $R_\varepsilon$, we consider the extension operator  \eqref{ext}.

\noindent The operator $\mathcal{A}_\varepsilon$ is defined as: $\mathcal{A}_\varepsilon f=u_f^\varepsilon$, where $f\in\mathcal{H}_\varepsilon$ and $u_f^\varepsilon\in\mathcal{V}_\varepsilon\subset\mathcal{H}_\varepsilon$ is the unique solution of the problem
	\begin{equation}
        \Phi_\ee(u_f^\ee , v )=
        \int_G fvdy\quad\forall v\in\mathcal{V}_\varepsilon.
        \label{ProblemaOperadorAeNeumann}
    \end{equation}
The above  existence and uniqueness of solution follows from the Lax-Milgram Lemma.

\noindent
Also it is simple to show that  $\mathcal{A}_\varepsilon$ is continuous, compact, self-adjoint, and positive. Indeed,  taking $v=u_f^\ee$ in   \eqref{ProblemaOperadorAeNeumann}, the Cauchy-Schwartz inequality leads us to
\begin{equation}
            \Phi_\varepsilon(u_f^\varepsilon,u_f^\varepsilon)^{\frac{1}{2}}\leq ||f||_{L^2(G)};
            \label{Desigualdad_FormaNeumann}
\end{equation}
        hence, to
        $$||\mathcal{A}_\varepsilon f||_{L^2(G)}=||u_f^\varepsilon||_{L^2(G)}\leq \Phi_\varepsilon(u_f^\varepsilon,u_f^\varepsilon)^{\frac{1}{2}}\leq ||f||_{L^2(G)},$$
        and, consequently,  the uniform bound  in \eqref{unif1} holds  true.

\noindent
For compactness, we use again the compact embedding of $H^1(G)$ into $L^2(G)$ and rewrite the proof in Section \ref{sec:ben3}, cf. \eqref{c1},  with minor modifications.

\noindent
The operator $\mathcal{A}_0$  is defined as: $\mathcal{A}_0 f=u_f$, where $f\in\mathcal{H}_0$ and $u_f\in\mathcal{V}_0$ is the unique solution of the problem
    \begin{equation}\label{lf}
        \Phi_0(u_f,v)=\langle f, v\rangle_{\mathcal{H}_0} \quad \forall v \in \mathcal{V}_0.
    \end{equation}
The above existence and uniqueness of the solution $u_f$ follows again from Lax-Milgram Lemma,  and the proof of continuity, compactness, self-adjointness and positiveness for operator  $\mathcal{A}_0$  follows as that for $\mathcal{A}_\varepsilon$.

\medskip
\noindent
Next, let us prove conditions $C_1-C_4$ in Lemma \ref{lemaole}.

\noindent
The verification of  $C_1$ and $C_2$ is performed by rewriting those in Section \ref{sec:ben3} with minor modifications.

\noindent
As regards  the proof of $C_3$, we proceed as follows.
Let us fix $f\in\mathcal{V}_0=H^1(l_0,l_1)$,  and   show \eqref{C3}  for  $\mathcal{H}_\varepsilon=L^2(G)$.
Indeed, considering $u_f^\varepsilon:=\mathcal{A}_\varepsilon R_\varepsilon f$ the unique solution of  \eqref{ProblemaOperadorAeNeumann} for $f\equiv R_\varepsilon f$, we  take $v=u_f^\varepsilon$ in  \eqref{ProblemaOperadorAeNeumann}
\begin{equation}
        \Phi_\varepsilon(u_f^\varepsilon,u_f^\varepsilon)=\langle R_\varepsilon f,u_f^\varepsilon\rangle_{\mathcal{H}_\varepsilon},
        \label{PhieComprobacionC3Neumann}
    \end{equation}
and obtain \eqref{Desigualdad_FormaNeumann} which implies the existence of  converging  subsequences of  $\{u_f^\varepsilon\}_{\ee} $ in $H^1$-weak, as $\ee\to 0$, and
\begin{equation}
||\partial_{y_1}u_f^\varepsilon||_{L^2(G)}^2+\varepsilon^{-2}||\partial_{y_2}u_f^\varepsilon||_{L^2(G)}^2+
    \varepsilon^{-2}||\partial_{y_3}u_f^\varepsilon||_{L^2(G)}^2\leq ||f||_{L^2(G)}^2.\label{convergenciaderivaday3Neumann}
\end{equation}
Hence, the bounds \eqref{convergenciaderivaday2} hold true and
the limit possible  functions of the subsequences $u_f^\varepsilon$ do not depend on the variables $y_2$ and $y_3$.

\noindent
As for the mixed boundary problem, we have that the sequence $u_f^\varepsilon$ is bounded and, therefore, from each sequence there exists a subsequence $\varepsilon'$ such that $u_f^{\varepsilon'}\to u$ in $H^1(G)$-weak and strong in $L^2(G)$ for some $u\in H^1(G)$, as $\ee'\to 0$. Let us  identify $u$ with $u_f^0=R_\varepsilon\mathcal{A}_0f$. To this end, we take functions $v_\varepsilon\in H^1(G)$ of the form
    $$v_\varepsilon(y_1,y_2,y_3):=v(y_1)\quad\text{ with }v\in\mathscr{C}^\infty(l_0,l_1).$$
    We have that
    $$\Phi_\varepsilon(u_f^\varepsilon,v_\varepsilon):=\int_G \partial_{y_1}u_f^\varepsilon\partial_{y_1}vdy+\int_G u_f^\varepsilon v dy =\langle R_\varepsilon f,v_\varepsilon\rangle_{L^2(G)}, $$
    and taking limits in the above equality, as $\ee'\to 0,$
    \begin{equation*}
        \int_G \partial_{y_1}u\partial_{y_1}vdy+\int_G u v dy =\int_{l_0}^{l_1}|D_{x_1}|\partial_{x_1}u\partial_{x_1}vdx_1+\int_{l_0}^{l_1}|D_{x_1}|uvdx_1 = \langle f,v\rangle_{\mathcal{H}_0}.
    \end{equation*}
    That is to say, $$\Phi_0(u,v)=\langle f,v\rangle_{\mathcal{H}_0}, $$
  while a density argument allows us to conclude that
    $$\Phi_0(u,v)=\langle f,v\rangle_{\mathcal{H}_0}\quad\forall v\in H^1(l_0,l_1).$$
 Therefore, $u$    satisfies \eqref{lf} and by the uniqueness of solution we have that the whole sequence converges towards  $u=\mathcal{A}_0f$. That is, $u_f^\varepsilon\to\mathcal{A}_0 f$ in $L^2(G)$, as $\ee\to 0$, and by definition it coincides with  $R_\varepsilon \mathcal{A}_0 f $. This concludes the proof of $C_3$.

\noindent
Finally, we verify condition $C_4$. Let $\{f_\varepsilon\}_{\varepsilon>0}$ be a bounded sequence in $\mathcal{H}_\varepsilon=L^2(G)$, that is, $\sup_{\varepsilon}||f_\varepsilon||_{L^2(G)}<\infty$. We consider the sequence $u_\varepsilon:=\mathcal{A}_\varepsilon f_\varepsilon$, which belongs to $H^1(G)$; that is, $u_\varepsilon$ is the solution of
    \begin{equation*}
        \Phi_\varepsilon(u_\varepsilon,v)=\langle  f_\varepsilon,v\rangle_{\mathcal{H}_\varepsilon}\quad\forall v\in H^1(G).
    \end{equation*}
Since the operators $\mathcal{A}_\varepsilon$ are uniformly bounded and the sequence $\{f_\varepsilon\}_{\varepsilon>0}$ is bounded, then the sequence $\{u_\varepsilon\}_{\varepsilon>0}$ defined above is bounded in $L^2(G)$, that is, $||u_\varepsilon||_{\mathcal{H}_\varepsilon}\leq C$ for some constant $C>0$. Then, rewriting the arguments for condition $C_3$ with minor modifications (see \eqref{PhieComprobacionC3Neumann}--\eqref{convergenciaderivaday3Neumann}), we obtain the bound $\Phi_\varepsilon(u_\varepsilon,u_\varepsilon)\leq K$ and the possible limits of $\{u_\varepsilon\}_{\varepsilon>0}$ are independent of $y_2$ and $y_3$.

\noindent
In this way, we can extract a subsequence $\{u_{\varepsilon'}\}_{\varepsilon'>0}$ converging strongly in $L^2(G)$ and weakly in $H^1(G)$:
    $$u_{\varepsilon'}\to u\quad \text{ in }L^2(G)\text{ for some }u\in H^1(G).$$
    We can identify $u$ with $R_{\varepsilon'} w_0$ for some $w_0\in\mathcal{V}_0=H^1(l_0,l_1)$, which concludes the verification of condition $C_4$.

    \noindent
Finally, it is worth noting that the eigenvalues $\{\tilde\mu_\varepsilon^n\}_{n=1}^\infty$ of operator $\mathcal{A}_\varepsilon$ are the inverses of the eigenvalues of problem \eqref{formulaciondebilproblemamixtoescaladoNeumann}, namely,
    \begin{equation}
        \tilde{\mu}_\varepsilon^n=\frac{1}{\tilde{\lambda}_\varepsilon^n}
        =\frac{1}{\lambda_\varepsilon^n+1} \quad n=1,2,\cdots.
        \label{i1N}
    \end{equation}
Similarly, the eigenvalues $\{\tilde\mu_0^n\}_{n=1}^\infty$ of $\mathcal{A}_0$ are the inverses of the eigenvalues of \eqref{formulaciondebilproblemalimiteNeumann}
\begin{equation}
	\tilde{\mu}_0^n=\frac{1}{\tilde{\lambda}_0^n}=\frac{1}{\lambda_0^n+1} \quad n=1,2,\cdots .
\label{i2N}
\end{equation}

\medskip
\noindent
Gathering all the results above, we apply Lemma \ref{lemaole} and, considering \eqref{i1N} and \eqref{i2N},  we have proved the convergence of the spectrum of \eqref{formulacionproblemaNeumann} towards that of \eqref{limitEDP},\eqref{limitN} with conservation of the multiplicity as well as the convergence of the corresponding eigenfunctions  as we state in the following theorem.

\begin{theorem}\label{TheoremN}
Suppose that the area of the cross sections of the domain $G$, $|D_{x_1}|$, is a piecewise continuous function for $x\in(l_0,l_1)$ satisfying condition \eqref{acotacionarea}. Let $\{\lambda_\varepsilon^n\}_{n=1}^\infty$ be the eigenvalues of problem \eqref{formulaciondebilproblemaNeumann} and $\{U_\varepsilon^n\}_{n=1}^\infty$ eigenfunctions forming an orthonormal basis in $L^2(G)$. Then, for any $n\in\mathbb{N}$, we have that
\begin{equation}\label{conservation2}
\lambda_\varepsilon^n\to\lambda_0^n\quad\text{ when }\varepsilon\to 0,
\end{equation}
where $\lambda_0^n$ is the $n$-th eigenvalue of the weak formulation of \eqref{limitEDP},\eqref{limitN} (cf. \eqref{formulaciondebilproblemalimiteNeumann}). Moreover, for each sequence $\varepsilon$, we can extract a subsequence (still denoted by $\varepsilon$) such that
\begin{equation}\label{conservation2fp}U_\varepsilon^n\to U_*^n\quad\text{in }H^1(G)\text{-weak when }\varepsilon\to 0,\end{equation}
where $U_*^n$ is an eigenfunction associated with $\lambda_0^n$ and such that $\{U_*^n\}_{n=1}^\infty$ form a basis of $H^1(l_0,l_1)$ and $L^2(l_0,l_1)$.
\\
In addition, for each $U_0$ eigenfunction of the limit problem associated with $\lambda_0^{k+1}$ eigenvalue with multiplicity $m$, $\lambda_0^{k}<\lambda_0^{k+1}=\lambda_0^{k+2}=\dots=\lambda_0^{k+m}<\lambda_0^{k+m+1}$, with $||U_0||_{\mathcal{H}_0}=1$, there exists a linear combination $\bar{U}_\varepsilon$ of the eigenfunctions $U_\varepsilon^{k+1},\dots,U_\varepsilon^{k+m}$ of problem \eqref{formulacionproblemaNeumann} such that
\begin{equation*}
        ||\bar{U}_\varepsilon-U_0||_{L^2(G)}\to 0\quad\text{ when }\varepsilon\to 0.
\end{equation*}
\end{theorem}

\section{Concluding Remarks}
\label{Darticulo}

In this research, we  tackle the asymptotics of singulary perturbed spectral problems for rod structures which remained as open problems in the literature of the Applied Mathematics. We also provide    a general mathematical framework that can be applied to many  other related models.

\medskip
\noindent
In particular, we  consider spectral problems for the Laplace operator in  3D rod structures with a  small cross section of diameter $O(\varepsilon)$, $\varepsilon$ being a positive parameter (see Figures \ref{Fig_domain}-\ref{figura_boveda}). The boundary conditions are Dirichlet  (Neumann, respectively)  on  the bases of this structure and Neumann on the lateral boundary. See problems  \eqref{formulacionuno_a},\eqref{formulacionuno_b},\eqref{fD} and  \eqref{formulacionuno_a},\eqref{formulacionuno_b},\eqref{fN}, respectively.  As   $\varepsilon\to 0$, we show the convergence of the eigenvalues with conservation of the multiplicity towards those of a 1D spectral  model  with Dirichlet (Neumann, respectively) conditions;  see \eqref{limitEDP},\eqref{limitD} and \eqref{limitEDP},\eqref{limitN} respectively.    We deal with the low frequencies and  the approach to eigenfunctions  in the suitable Sobolev spaces is also shown, providing information on the longitudinal oscillations which take place at a low-frequency level, as $\ee\to 0$. As a sample,  we provide some graphics of approximations of eigenfunctions obtained by numerical methods (see Figures \ref{figura_mixed}-\ref{figura_mixed2}). However,  having observed numerical instabilities as $\varepsilon$ becomes smaller (see \cite{articulo} for further details), the asymptotic analysis that we perform here becomes necessary. Let us refer to the convergence \eqref{conservation},\eqref{conservationfp} for the mixed boundary value problem   and \eqref{conservation2},\eqref{conservation2fp} for the Neumann one.

\medskip
\noindent
The limit 1D problems may arise in diffusion or vibration models of nonhomogeneous media with different physical characteristics and they take into account the geometry of the 3D domains.  The interest of these 3D models from the dynamical viewpoint is also evident (cf. Section \ref{sec:ben1}). However, we also address convergence of solutions for associated stationary models, which were open problems in the literature: see problems \eqref{ProblemaOperadorAe} and \eqref{Phi0ComprobacionC3},     and \eqref{ProblemaOperadorAeNeumann}  and \eqref{lf} for the different associated bilinear forms.

\medskip
\noindent
Among the techniques here used, we mention dimension reduction procedures and techniques from the spectral perturbation theory; in particular, we are led to a fixed domain of reference by introducing  the stretching variables (cf. \eqref{dominioG}-\eqref{acotacionarea} and \eqref{cambioGe-G}).

\medskip
\noindent
Finally, let us note that the explicit and numerical computations in \cite{articulo} have also been  performed for the Dirichlet problem enlightening that the oscillations for the eigenfunctions associated to the low frequencies hold in the transverse directions: the study  of both  low and high frequencies remaining as an open problem to be considered by the authors. In contrast, the study of the eigenvalues when the Neumann condition is imposed on the lateral boundary and one of the bases,    and Dirichlet condition on the another one,   seems to follow the structure in this paper, the limit problem being \eqref{limitEDP} with mixed boundary conditions  on the end points. For brevity we avoid   this study here; cf. \cite{CardoneDranteNazarov2010} for different boundary conditions.

\medskip
\noindent
{\bf Acknowledgments:}
This work has been partially supported by Grant PID2022-137694NB-I00 funded by MICIU/AEI/ 10.13039/501100011033 and by ERDF/EU.

\end{document}